\documentclass[11pt]{article}

\usepackage[T1]{fontenc}
\usepackage[latin1]{inputenc}
\usepackage{mathtools,amsthm,amssymb}
\usepackage{mathabx}

\usepackage{authblk}

\usepackage{xcolor}
\usepackage{graphicx,caption,subcaption}
\graphicspath{{figs/}}

\usepackage{hyperref}

\usepackage[margin=1in]{geometry}

\newtheorem{prop}{Proposition}[section]
\newtheorem{cor}[prop]{Corollary}

\newtheorem{lem}[prop]{Lemma}
\newtheorem{theo}[prop]{Theorem}

\numberwithin{equation}{section}

\newcommand*{\longhookrightarrow}{\ensuremath{\lhook\joinrel\relbar\joinrel\rightarrow}}

\newcommand{\vertiii}[1]{{\left\vert\kern-0.25ex\left\vert\kern-0.25ex\left\vert #1 
    \right\vert\kern-0.25ex\right\vert\kern-0.25ex\right\vert}}

\def\tr{\mbox{\rm Tr}}

\newcommand{\EE}{\mathbb{E}}

\newcommand{\HH}{\mathbb{H}}

\newcommand{\LL}{\mathbb{L}}
\newcommand{\MM}{\mathbb{M}}
\newcommand{\NN}{\mathbb{N}}
\newcommand{\PP}{\mathbb{P}}

\newcommand{\RR}{\mathbb{R}}

\newcommand{\XX}{\mathbb{X}}

\newcommand{\Ca}{ {\cal C }}

\newcommand{\La}{ {\cal L }}

\newcommand{\Ea}{ {\cal E }}
\newcommand{\Sa}{ {\cal S }}

\newcommand{\Va}{ {\cal V }}
\newcommand{\Ua}{ {\cal U }}
\newcommand{\Fa}{ {\cal F }}

\newcommand{\Ia}{ {\cal I }}
\newcommand{\Xa}{ {\cal X }}

\newcommand{\Ta}{ {\cal T}}
\newcommand{\Ha}{ {\cal H }}

\newcommand{\Pa}{ {\cal P }}

\newcommand{\Ya}{ {\cal Y }}
\newcommand{\Wa}{ {\cal W }}

\def \PP{\mathbb{P}}
\def \RR{\mathbb{R}}
\def \SS{\mathbb{S}}
\def \EE{\mathbb{E}}

\def \LL{\mathbb{L}}

\newcommand{\point}{\mbox{\LARGE .}}

\newcommand{\cqfd}{\hfill\blbx \medskip}
\def\blbx{\hbox{\vrule height 5pt width 5pt depth 0pt}\medskip}

\makeatletter
\DeclareRobustCommand\frownotimes{\mathbin{\mathpalette\frown@otimes\relax}}
\newcommand{\frown@otimes}[2]{
  \vbox{
    \ialign{##\cr
      \hidewidth$\m@th#1{}_\frown$\kern-\scriptspace\hidewidth\cr
      \noalign{\nointerlineskip\kern-.1pt}
      $\m@th#1\otimes$\cr
    }
  }
}
\makeatother

\begin{document}

\title{A perturbation analysis of stochastic matrix Riccati diffusions}

\author[$1$]{Adrian N. Bishop}
\author[$2$]{Pierre Del Moral}
\author[$3$]{Ang\`ele Niclas}
\affil[$1$]{{\small University of Technology Sydney (UTS) and CSIRO}}
\affil[$2$]{{\small INRIA, Bordeaux Research Center, France}}
\affil[$3$]{{\small \'Ecole Normale Sup\'erieure de Lyon}}

\date{}

\maketitle

\begin{abstract}
Matrix differential Riccati equations are central in filtering and optimal control theory. The purpose of this article is to develop a perturbation theory for a class of stochastic matrix Riccati diffusions. Diffusions of this type arise, for example, in the analysis of ensemble Kalman-Bucy filters since they describe the flow of certain sample covariance estimates. In this context, the random perturbations come from the fluctuations of a mean field particle interpretation of a class of nonlinear diffusions equipped with an interacting sample covariance matrix functional. The main purpose of this article is to derive non-asymptotic Taylor-type expansions of stochastic matrix Riccati flows with respect to some perturbation parameter. These expansions rely on an original combination of stochastic differential analysis and nonlinear semigroup techniques on matrix spaces. The results here quantify the fluctuation of the stochastic flow around the limiting deterministic Riccati equation, at any order. The convergence of the interacting sample covariance matrices to the deterministic Riccati flow is proven as the number of particles tends to infinity. Also presented are refined moment estimates and sharp bias and variance estimates. These expansions are also used to deduce a functional central limit theorem at the level of the diffusion process in matrix spaces.\\
\\
{\em R\'esum\'e:}
Les \'equations de Riccati matricielles jouent un r\^ole important dans la th\'eorie du filtrage et du contr\^ole optimal. Cet article pr\'esente une th\'eorie des perturbations d'une classe d'\'equations  de Riccati matricielles stochastiques. Ces mod\`eles probabilistes sont d'un usage courant dans la th\'eorie des filtres de Kalman d'Ensemble. Ils repr\'esentent dans ce contexte l'\'evolution des matrices de covariance empiriques associ\'ees \`a un ensemble de diffusions en interaction. Les perturbations al\'eatoires r\'esultent des fluctuations stochastiques d'un syst\`eme de particules de type champ moyen interagissant avec la mesure empirique du syst\`eme. Nous pr\'esentons dans cet article une formule de Taylor non asymptotique pour des flots stochastiques de diffusion de Riccati matricelles par rapport \`a un param\`etre de fluctuation. Ces d\'eveloppements sont fond\'es sur un nouveau calcul diff\'erentiel stochastique et une analyse fine de semigroupes non lin\'eaires dans des espaces de matrices. Ces r\'esultats permettent de quantifier avec pr\'ecision les fluctuations des flots de matrices stochastiques autour des syst\`emes limites à tout ordre. Nous illustrons ces r\'esultats avec une preuve de la convergence des matrices empiriques de filtres de Kalman d'Ensemble vers la solution d'\'equations de Riccati d\'eterministes lorsque le nombre de particules tends vers l'infini. Nous pr\'esentons dans ce cadre des estimations fines des biais et des variances, ainsi qu'un theor\`eme de la limite centrale fonctionnel au niveau du processus matriciel.
\end{abstract}

\section{Introduction}\label{sec-intro}

Matrix Riccati equations play a central role in stochastic filtering and optimal control theory. These quadratic differential equations are used to design optimal Kalman filters and optimal controllers in dual quadratic cost and linear system regulation problems. 

This article presents a perturbation and fluctuation analysis for a class of matrix diffusions combining a Riccati drift functional with a diffusive martingale with a cubic-type predictable angle bracket. This class of stochastic model is defined in terms of $\Wa_t$: a $(r\times r)$-matrix with independent Brownian entries with $r\geq 1$.  We associate with some diagonalizable
$(r\times r)$-matrix $A$, and some positive definite matrices $R,S> 0$, the Riccati drift function $\Lambda$ defined by
\begin{equation}\label{def-Riccati-drift}
	\Lambda(Q) :=(A-QS)Q+Q(A-QS)^{\prime}+\Sigma(Q)\quad \mbox{\rm with}\quad \Sigma(Q):=R+QSQ
\end{equation}
The Riccati diffusions discussed in this article are then defined by the stochastic model
\begin{equation}\label{f21}
	dQ_t=\Lambda(Q_t)\,dt+\epsilon\,dM_t\quad\mbox{\rm with}\quad dM_t:=\left[\varphi(Q_t)~d\Wa_t~\Sigma_{\varphi}\left(Q_t\right)\right]_{ sym}
\end{equation}
In the above display, $B^{\prime}$ stands for the transposition of $B$,  $B_{ sym}=(B+B^{\prime})/2$ the symmetric part of an $(r\times r)$-matrix, 
$\varphi$ is square root function $\varphi(Q):=Q^{1/2}$ and  $\Sigma_{\varphi}:=\varphi\circ\Sigma$.

We let $\phi^{\epsilon}_{t}(Q_0)=Q_{t}$ be the stochastic flow associated with the solution $Q_t$ of the equation (\ref{f21}) starting at $Q_0$. 
To clarify the presentation, we also write  $\phi_{t}$  instead of $\phi_{t}^{0}$, i.e. the semigroup  associated with the Riccati equation (\ref{f21}) when $\epsilon=0$.
 
The diffusion term relies on some parameter $\epsilon\in [0,1]$ which reflects the variance of the perturbations.  
The fluctuation parameter $\epsilon$ is chosen so that $2\epsilon^2r<1$ to ensure the existence and the boundedness properties of these quadratic diffusion processes (cf. (\ref{trace-estimates-unif}) in theorem~\ref{th1-estimates-Ln}). When $\epsilon=0$, the semigroup associated with these models resumes to the conventional unperturbed matrix Riccati equation. The fluctuation of the diffusion around the unperturbed Riccati equation can be quantified in terms of 
$\LL_n$-mean error estimates, as soon as $5n\epsilon^2r<1$ (cf. (\ref{Frob-Q-Q-estimates}) in theorem~\ref{th1-estimates-Ln}).

Whenever $S=0$ the evolution equation (\ref{f21}) resumes to the Wishart process. The case $A=0$ is also known as the squared Brownian motion (a.k.a. matrix square Bessel processes). 

Related diffusions in symmetric matrix spaces arise in a variety of application domains. For example, backward stochastic matrix Riccati equations arise in linear quadratic optimal control problems with random coefficients~\cite{bismut1976,Kohlmann2003,hu-zhou}. The class of stochastic processes in matrix spaces discussed in this article also encapsulates conventional Wishart processes arising in mathematical physics~\cite{katori1,katori2}, multivariate statistics~\cite{barndorff,graczyk}, econometrics and financial 
mathematics~\cite{cuchiero,gourieroux}. In this context, the parameter $\epsilon^2$ captures the amplitude of the fluctuations of the volatility process. 

This article may be motivated by applications in signal processing and in filtering and data assimilation in high-dimensional inference problems. In this context, a stochastic matrix Riccati equation (like (\ref{f21})) represents the evolution of the sample covariance matrices associated with the ensemble Kalman-Bucy filter~\cite{delmoral16a} (abbreviated {\tt EnKF}). In this case, the parameter $\epsilon^2$ is inversely proportional to the number of particles associated with these filters. This class of {\tt EnKF} should be interpreted as a mean-field particle approximation of a nonlinear McKean-Vlasov-type diffusion; see the latter application of our theory to the {\tt EnKF}, and also \cite{delmoral16a,dbkr-2017,Bishop/DelMoral/multiDimRicc}. These general probabilistic models were introduced by H.P. McKean~\cite{mckean}; see also \cite{sznitman,meleard2,dm-crc2013} for a detailed discussion and applications of these general models as well as the more recent article~\cite{dm-a} dedicated to the long time behavior of interacting diffusions. Since its introduction in the early 1990s, the {\tt EnKF} has been widely studied and applied for numerically solving forecasting and data assimilation problems \cite{evensen03,houte,kalnay}. One interesting numerical aspect of the {\tt EnKF} is that they can be adapted to work well in high-dimensional, nonlinear, and small-noise scenarios. 

Under appropriate observability and controllability conditions, we know that the true (classical) Kalman filter can track (unstable) linear noisy signals uniformly w.r.t. the time horizon (cf.~\cite{ap-2016,bd-CARE} and references therein), and that the conditional error covariance satisfies a stable differential matrix Riccati equation. This deterministic Riccati equation and its limiting behaviour is actually central in the stability analysis of the classical Kalman filter \cite{ap-2016}. One currently open research stream concerns the {\tt EnKF} and its ability to also track any unstable modes of a noisy signal in any dimension (under appropriate regularity conditions) with noisy observations.

In the context of linear-Gaussian filtering problems, the stability and convergence properties of the {\tt EnKF} rely heavily on the fluctuations of a particular stochastic matrix Riccati diffusion \cite{delmoral16a,Bishop/DelMoral/multiDimRicc}. This equation models the flow of the {\tt EnKF} sample covariance matrix (analogously to how a related differential matrix Riccati equation models the flow of the true covariance matrix in classical Kalman filtering \cite{ap-2016}). 

The present article focuses on this stochastic matrix Riccati diffusion and its perturbation and fluctuation properties. These results are of broad mathematical interest on their own, since the Riccati diffusion considered is a rather general quadratic matrix-valued stochastic differential equation, see also \cite{Bishop/DelMoral/multiDimRicc}. These results are also of interest under the {\tt EnKF} banner, as they rigorously characterise some behaviour of the {\tt EnKF} sample covariance flow. 

We present sharp and non-asymptotic expansions of matrix moments of the matrix Riccati diffusion with respect to the parameter $\epsilon$, stripped of all analytical superstructure, and probabilistic irrelevancies. These results can basically be stated as follows
\begin{equation}\label{intro-main-Taylor}
\phi^{\epsilon}_t=\phi_t+\sum_{1\leq k<n}~\frac{\epsilon^k}{k!}~\partial^k\phi_t+\overline{\partial}^{\,n}\!\phi^{\epsilon}_t
\end{equation}
for any $n\geq 1$ and some stochastic flow $\partial^k\phi_t$ {\em whose values don't depend on the fluctuation parameter $\epsilon$} and some remainder stochastic term $\overline{\partial}^{\,n}\phi^{\epsilon}_t$ of order $\epsilon^n$. We also provide uniform estimates of the stochastic flow $\partial^k\phi_t$ w.r.t. the time horizon {\em even when the matrix $A$ is unstable} (see (\ref{estimate-partial-k}) in theorem~\ref{theo-intro-phi-MM-Taylor}, the estimate (\ref{monotone-prop}) and section~\ref{stability-riccati-square-root}).

These estimates are stronger than
the conventional functional central limit theorems for stochastic processes. For example, these results clearly imply the almost sure central limit theorem 
$$
\epsilon^{-1}\left[\phi_t^{\epsilon}-\phi_t\right]~\longrightarrow_{\epsilon\rightarrow 0}~\partial\phi_t
$$

These matrix moment expansions rely on an original combination of stochastic differential analysis with nonlinear semigroup techniques on matrix spaces. The Taylor-type matrix representation discussed in this article allows one to quantify the fluctuation at any order of the stochastic flow around the limiting deterministic Riccati equation. Moreover, they also provide refined moment estimates as well as sharp bias and variance estimates. For instance, using (\ref{intro-first-bias}) and theorem~\ref{th2} we have the bias estimate
$$
\phi_t(Q)+\frac{\epsilon^2}{2}~\EE\left(\partial^2\phi_t(Q)\right)+\mbox{\rm O}(\epsilon^4)= \EE\left(\phi^{\epsilon}_t(Q)\right)~\leq~  \phi_t(Q)
$$
An explicit description of the bias matrix flow $t\mapsto\EE\left(\partial^2\phi_t(Q)\right)$ is provided in theorem~\ref{th2}.
Uniform estimates w.r.t. the time parameter of the derivative processes are also provided; e.g. at every time step, or at the level of the path. Last but not least, we also extend these expansions to stochastic flows starting from random fluctuation matrices. Combined with the Fa\`a  di Bruno's formula they can be used to deduce any matrix moments and any matrix-moment of smooth functionals of the stochastic Riccati flows.

{The article is organized into four parts:} 

\noindent The first part in Section \ref{sec-intro} is dedicated to the description of stochastic matrix Riccati diffusions and a statement of the main results (Section \ref{main-results-sec}). The primary study in this work is a detailed fluctuation analysis of these matrix Riccati diffusions. In Section \ref{sec-EKF} we illustrate the main results and their significance in the context of the {\tt EnKF}. In Section~\ref{stability-riccati-square-root}, we briefly discuss the stability and long-time behaviour of these Riccati diffusions. We describe the invariant measure of one-dimensional stochastic Riccati diffusions, and we relate some moment explosion properties exhibited by these diffusions with the heavy-tailed structure of this stationary measure. The stability exposition in Section~\ref{stability-riccati-square-root} is further detailed in depth in \cite{dbkr-2017,Bishop/DelMoral/multiDimRicc}.

The second part of the article in Section \ref{smoothnes-section} is concerned with the smoothness properties of the semigroups associated with stochastic matrix Riccati diffusions. We also provide a brief review on matrix analysis and matrix functional differential techniques, including a Fa\`a  di Bruno's formula with remainder and Taylor expansions of square root functionals.

The third part in Section \ref{stoch-matrix-calc-section} is dedicated to stochastic matrix integration and matrix valued martingales. We describe general formulae for computing the predictable angle brackets of matrix valued martingales in terms tensor and symmetric tensor products. We also provide an Ito formula for matrix functionals of stochastic Riccati equations and a series of martingale continuity theorem.

The last part in Section \ref{proof-section} is concerned with the technical proof of the four main theorems presented in this article. The proofs of some technical lemmas are provided later in Appendix \ref{appendix-section}. 

\subsection{Some basic notation}

We introduce some common notation used throughout this article, and in particular in the statement of our main results given subsequently. Later in Section \ref{smoothnes-section} we provide further notation related to the developments in that section concerning the smoothness properties of related semigroups and specific Taylor-type expansions. In Section \ref{stoch-matrix-calc-section} we detail specific notation relevant to the developments of that section which is devoted to stochastic matrix calculus and matrix martingale theory. Rarely, and locally within (sub-)sections, some symbols may be reused where there is no risk of confusion.

We denote by $\Ta_{r,r^{\prime}}$ the space of $(r\times r^{\prime})$-matrices with real entries and some $r,r^{\prime}$. We also let $\Sa_r\subset \Ta_r:=\Ta_{r,r}$ denote the closed subset of symmetric matrices, $\Sa_r^0\subset\Sa_r$ the subset of positive semi-definite matrices, and $\Sa_r^+\subset \Sa_r^0$ the open subset of positive definite matrices. Given $Q\in \Sa_r^0-\Sa_r^+$ we denote by $Q^{1/2}$ a (non-unique) but symmetric square root of $Q$ (given by a Cholesky decomposition). When $Q\in\Sa_r^+$ we always choose the principal (unique) symmetric square root. With a slight abuse of notation we denote by $I$ the $(r\times r)$ identity matrix, for any $r\geq 1$. In the further development of this article $\Vert\point\Vert$ denotes the spectral norm $\|\cdot\|_2$ or the Frobenius norm $\|\cdot\|_F$ on the space of matrices. We let $\rho(A)$ denote the logarithmic norm of the matrix $A$. Finally, for any square matrix $A$ we define a particular matrix operation notation $\{\cdot\}$ by
\begin{equation}\label{def-A-par}
 \{A\}=2^{-1}\left(A+\tr(A)\,I\right)
\end{equation}
where $\tr(\cdot)$ denotes the trace operator. 

We let  $\Ca\left([0,t],\Ta_{r}\right)$ be the space of continuous stochastic flows $s\in [0,t]\mapsto A_s$ from $[0,t]$ into the space $\Ta_{r}$ of $(r\times r)$-matrices. For any $m,n\geq 1$,  we equip  $\Ca\left([0,t],\Ta_{r}\right)^m$ with the norms defined for any flow $A_s:=(A_s^{(i)})_{1\leq i\leq m}$ of stochastic matrices by
$$
\vertiii{A}_{t,n}:=
 \sum_{1\leq i\leq m}\vertiii{A^{(i)}}_{t,n}\quad\mbox{\rm with}\quad\vertiii{A^{(i)}}_{t,n}=\EE\left[\Vert A^{(i)}\Vert_{t}^n\right]^{1/n}\quad\mbox{\rm and}\quad
\Vert A^{(i)}\Vert_{t}:=\sup_{s\in [0,t]}\Vert A^{(i)}_s\Vert
$$
For $m=1$ and for matrices $A_t$, we also consider the $\LL_n$-norm given by $\vertiii{A_t}_{n}=\EE\left[\Vert A_t\Vert^n\right]^{1/n}$. 

A random mapping $\epsilon\in [0,1]\mapsto A^{\epsilon}\in \Ca\left([0,t],\Ta_r\right)$ is said to be $n$-th differentiable as soon as for any $s\in[0,t]$ we have 
$$
 A^{\epsilon}_s=A_s+\sum_{1\leq k<n}~\frac{\epsilon^k}{k!}~\partial^kA_s+\overline{\partial}^{\,n}A^{\epsilon}_s
$$
for some stochastic flows processes $\partial^kA\in \Ca\left([0,t],\Ta_r\right)$ (whose values doesn't depends on $\epsilon$) and some 
remainder term $\overline{\partial}^{\,n}A^{\epsilon}_s$ such that 
$$
\forall m\geq 1\quad\exists \epsilon_m\in [0,1]\quad \mbox{\rm s.t.}\quad \forall \epsilon\in \left[0, \epsilon_m\right]\qquad
\vertiii{\overline{\partial}^{\,n}\!A^{\epsilon}}_{t,m}\leq a_{m,n}(t)~\epsilon^n
$$
for some finite parameters  $a_{m,n}(t)$. We also consider the collection of matrix-valued martingales 
$$
\MM^{\epsilon}_t(Q):=\int_0^t\left[
\varphi\left(\phi^{\epsilon}_{s}(Q)\right)\,d\Wa_s~\Sigma_{\varphi}\left(\phi^{\epsilon}_{s}(Q)\right)\right]_{ sym}
$$
To clarify the presentation, we also write  $\MM_t(Q)$ instead of  $\MM^{0}_t(Q)$.

We let $E_{s,t}(Q)$ be the exponential transition semigroup associated with the flow of matrices $\left[A-\phi_u(Q)S\right]$, $s\leq u\leq t$; a more explicit description is given in (\ref{estimate-Ea-t-bis}). When $s=0$, we write $E_{t}(Q)$ instead of $E_{0,t}(Q)$. With this notation, we have $E_{s,t}(Q)=E_{t}(Q)E_{s}(Q)^{-1}$. 

We also consider the matrix functionals $( \Gamma_t,\Pi_t,\Omega_t)$ from $\Sa_r^0$ into $\Sa_r$ defined by
$$
\begin{array}{rcl}
\displaystyle\Omega_t(Q)&:=&E_t(Q)\Pi_t(Q)E_t(Q)^{\prime} \\
\displaystyle\mbox{\rm with}\quad\displaystyle
\Pi_{t}(Q)&:=&Q\,\left\{\Gamma_{t}(Q)~ \Sigma(Q)\right\}+ \Sigma(Q)\,\left\{\Gamma_{t}(Q)~Q\right\} \quad\mbox{\rm and}\quad \displaystyle	\Gamma_t(Q):= \int_0^t~E_u(Q)^{\prime}\,E_u(Q)\,du 
\end{array}
$$

The article discusses several stochastic Taylor-type expansions and diffusion processes in matrix spaces. The stochastic analysis of these combines several combinatorial and algebraic sophisticated tools with Burkholder-Davis-Gundy inequalities and related $\LL_n$-error norms between stochastic process. We will not track the dependency on the parameter $n\geq 1$. We have chosen to focus on the dependencies w.r.t. the time horizon.

In the further development of the article  for any time horizon $t\geq 0$ we  set
$$
e(t):=c\exp{(k\rho(A) t)}~\leq~ e_+(t):=c~(1\vee \exp{(k\rho(A) t)}) ~\leq~ \overline{e}(t):=(1+t^l)~e_+(t)
$$
When $\rho(A)<0$ the function $e(t)$ tends to $0$ exponentially fast, while $e_+(t)$ is uniformly bounded and $\overline{e}(t)$ has polynomial growth. If $\rho(A)>0$, then $e(t)$ grows exponentially fast. Note that $\rho(A)<0$ is a (type of strong) stability condition on $A$, and is stronger than we suspect necessary (in general) for time-uniform boundedness in our main results below. For example, in one-dimension, time-uniform versions of theorems \ref{th1-estimates-Ln} and \ref{th2} hold under (much) weaker conditions; see \cite{dbkr-2017}. Less ``detailed'' fluctuation results in the matrix-valued setting are explored in \cite{Bishop/DelMoral/multiDimRicc} without this stability condition on $A$. In this article, our main focus is on a very complete (e.g. ``at any order'', and on the path space) perturbation and fluctuation analysis of the matrix valued diffusion (\ref{f21}), and we do not pursue further relaxations on the model here.

For any matrix  $Q\in \Sa^+_r$ we also consider the following parameters
$$
\vert Q\vert:=\vert Q\vert_-+\vert Q\vert_+\quad \mbox{\rm with}\quad \vert Q\vert_-:=c_1~(1+\Vert Q^{-1}\Vert^{n_1})
\quad \mbox{\rm and}\quad
 \vert Q\vert_+:=c_2~(1+\Vert Q\Vert^{n_2})
$$
In the above display formulae $n_1,n_2\geq 0$ and $c_1,c_2$ are some finite constants whose values may vary from line to line. 
Finally, we use the sign $\longhookrightarrow_{\epsilon\rightarrow 0}$ to denote the weak convergence of processes when 
$\epsilon$ tends to $0$. We also use the notation $a\vee b=\max{(a,b)}$ and  $a\wedge b=\min{(a,b)}$ for any $a,b\in\RR$, and the conventions
$\sum_{\emptyset}=0$ and $\prod_{\emptyset}=1$.

\subsection{Statement of the main results} \label{main-results-sec}
 
 As noted, this article is concerned primarily with the detailed perturbation and fluctuation properties of the matrix Riccati diffusion (\ref{f21}). We point to \cite{dbkr-2017,Bishop/DelMoral/multiDimRicc} for some related fluctuation results, and more particularly to some stability and contraction results on the Riccati diffusion.
 
Whenever $S>0$, up to a change of basis, there is no loss of generality to assume that $S=I$. More precisely the stochastic matrices $\overline{Q}_t:=S^{1/2}Q_tS^{1/2}$ satisfy the same equation as (\ref{f21}) when we replace $(A,R,S)$ by the matrices 
\begin{equation}\label{ref-overline-A}
	(\overline{A},\overline{R},\overline{S}):=(S^{1/2}AS^{-1/2},S^{1/2}RS^{1/2},I)
\end{equation}
The invariance property of the drift follows elementary algebraic manipulations. The analysis of the martingale part is a direct consequence of the formula (\ref{coupligng-eq-W}) presented later. In the further development of the article, unless otherwise stated, we assume that $S=I$. In this case, and with $R>0$ we may take the standard controllability and observability conditions as holding; see \cite{Antsaklis,ap-2016} for details on controllability and observability in control and filtering theory. 

 The first main result concerns a uniform bound on any differential moment given in terms of the minimum eigenvalue of the Riccati semigroup and the log-norm of $A$.
 
 \begin{theo}\label{theo-intro-phi-MM-Taylor}
 For any $Q\in \Sa_r^+$, the stochastic flows $\epsilon\mapsto \phi^{\epsilon}_t(Q)$
 and $\epsilon\mapsto \MM^{\epsilon}_t(Q)$ are smooth. For any $m,n\geq 1$, $t\geq 0$ and  
 for any $\delta>0$ we have the uniform estimates
\begin{equation}
\begin{array}{l}
\displaystyle\vertiii{\partial^n\phi_t(Q)}_{m}\vee\left[(1+t)^{-n/2}\left(\vertiii{\partial^n\phi(Q)}_{t,m}\vee \vertiii{\partial^{n-1}\MM}_{t,m}\right)\right] 
\leq ~c~\vert Q\vert_-~\exp{\left(\delta\Vert Q\Vert\right)} \label{estimate-partial-k}
\end{array}
\end{equation}
for some constant $c$ whose values only depend on  $(\delta,m,n)$.
In addition there exists some collection of parameters $\epsilon_{m,n}\in [0,1]$ such that for any $\epsilon\in [0, \epsilon_{m,n}]$ we have the remainder estimates
    \begin{equation}\label{intro-phi-MM-Taylor}
    \begin{array}{l}
\vertiii{\epsilon^{-n}\,\overline{\partial}^{\,n}\phi^{\epsilon}(Q)}_{t,m}\vee \vertiii{\epsilon^{-n+1}\,\overline{\partial}^{\,n-1}\MM^{\epsilon}(Q)}_{t,m}~\leq~ \overline{e}(t)~\vert Q\vert~
\end{array}
  \end{equation}
 \end{theo}

The proof of this theorem is provided in Section~\ref{proof-theo-intro-phi-MM-Taylor}. We have the  bias estimate
\begin{equation}\label{intro-first-bias}
	(\ref{intro-phi-MM-Taylor}) ~~~~\Longrightarrow  ~~~~  \Vert \EE\left(\phi^{\epsilon}_t(Q)\right)-\phi_t(Q)-2^{-1}~\epsilon^2~\EE\left(\partial^2\phi_t(Q)\right)
\Vert~\leq~ \epsilon^4~\overline{e}(t)~\vert Q\vert
\end{equation}
A detailed proof of (\ref{intro-first-bias}) is given toward the end of Section~\ref{proof-theo-intro-phi-MM-Taylor}, following the proof of theorem \ref{theo-intro-phi-MM-Taylor}.
      
Using the evolution equation (\ref{f21}), the computation of the relevant matrix differentials follows standard differential rules
on composition mappings and Fr\'echet derivation techniques. For instance, the first and second order derivatives of the flow $\phi^{\epsilon}_t(Q)$ are given by the formulae
\begin{eqnarray}
\displaystyle\partial\phi_t(Q)&=&
\int_0^t~\nabla\Lambda(\phi_s(Q))\cdot
\partial\phi_s(Q)
~ds+\MM_t(Q)\label{def-partial-phi}\\
\displaystyle\partial^{2}\phi_t(Q)&=&
\int_0^t
~\nabla^{2}\Lambda(\phi_s(Q))\cdot\left(\partial\phi_s(Q),\partial\phi_s(Q)\right)
~ds+
\int_0^t
~\nabla\Lambda(\phi_s(Q))\cdot\partial^{2}\phi_s(Q)
~ds+2~\partial\MM_t(Q)\nonumber
\end{eqnarray}
with the first order matrix-valued martingale
$$
\partial\MM_t(Q)=\int_0^t\left[(\varphi\circ\phi_s)(Q)~d\Wa_s~
\partial(\Sigma_{\varphi}\circ\phi_s)(Q)+\partial(\varphi\circ\phi_s)(Q)~d\Wa_s~
(\Sigma_{\varphi}\circ\phi_s)(Q)\right]_{ sym}\
$$
In the above display, $\nabla^k(\cdot)\cdot(\cdot,\ldots)$ stands for the $k$-th Fr\'echet derivative operator; see \cite{higham} and Section~\ref{sec-derivatives} for a more detailed discussion of these differentials. A sequential description of the derivatives at any order of the stochastic flow is given in Section~\ref{proof-theo-intro-phi-MM-Taylor} in the proof of theorem \ref{theo-intro-phi-MM-Taylor}.

The next theorem discusses some more refined first order estimates.
\begin{theo}\label{th1-estimates-Ln}
 For any $Q\in \Sa_r^0$ and for any $n\geq 1$ and $\epsilon\in [0,1]$ such that $2n\epsilon^2r<1$ we have 
\begin{eqnarray}
\vertiii{\phi^{\epsilon}(Q)}_{t,n}&\leq& (1+t)\,e_+(t)\,\vert Q\vert_+,\qquad{\vertiii{\phi^{\epsilon}_t(Q)}_{n}} \leq \vert Q\vert_+
\qquad~\label{trace-estimates-unif}
\end{eqnarray}
In addition, when $5n\epsilon^2r<1$  we also have the estimates
\begin{equation}\label{Frob-Q-Q-estimates}
\left[t^{-1/2}\vertiii {\phi^{\epsilon}(Q)-\phi(Q)}_{t,n}\right]\vee
     \vertiii {\phi_t^{\epsilon}(Q)-\phi_t(Q)}_{n}\leq \epsilon~e_+(t)
~\vert Q\vert_+
\end{equation}
\end{theo}

The proof of this theorem is provided in Section~\ref{th1-estimates-Ln-proof}. The estimates (\ref{Frob-Q-Q-estimates}) extend the ones presented in~\cite{delmoral16a} to path spaces and non necessarily stable matrices $A$. Also notice that the existence of the $n$-th moments (\ref{trace-estimates-unif}) of the stochastic flows $\phi^{\epsilon}_t(Q)$ requires a sufficiently small perturbation parameter. This property doesn't come from any technical overestimation, but from the heavy tailed properties of the stationary measures of the flows. A more thorough discussion on these properties and the possible moment explosion times of these stochastic flows is provided in Section~\ref{stability-riccati-square-root} in the context of one dimensional models.

The next theorem provides bias and second order estimates.
\begin{theo}\label{th2}
 For any $Q\in \Sa_r^0$ and any $n\geq 1$ there exists some $\epsilon_n\in [0,1]$ such that
for any time horizon $t\geq 0$ and any $\epsilon\in [0,\epsilon_n]$ we have
\begin{equation}\label{order-1-sobolev}
 \EE\left(\phi^{\epsilon}_t(Q)\right)\leq  \phi_t(Q)
\quad\mbox{and}\quad
\vertiii{\phi^{\epsilon}(Q)-\phi(Q)-\epsilon \,\partial\phi(Q)}_{t,n}\leq 
~\epsilon^{3/2} \overline{e}(t)~\vert Q\vert_+
\end{equation}
 In addition, we have the asymptotic bias formula
\begin{equation}\label{equivalent-bias-terms}
\EE\left[\partial^2\phi_t(Q)\right]=-\int_0^t \Omega_{t-s}\left(\phi_s(Q)\right)=-\int_0^t~E_{s,t}(Q)~\EE\left(\left[\partial\phi_s(Q)\right]^2\right)~E_{s,t}(Q)
~ds\leq 0
\end{equation}
\end{theo}

The proof of this theorem is provided in Section~\ref{th2-proof}. When the fluctuations come from some approximation scheme, such as in ensemble Kalman-Bucy filters, the l.h.s inequality in (\ref{order-1-sobolev}) shows that the stochastic flow under estimates the 
solution of the Riccati equation even when they start from the same initial matrix. That is, given any sufficient number of particles to ensure the diffusion is well-enough behaved, the sample covariance estimates computed in the {\tt EnKF} algorithm are always negatively biased.  

Now assume that $H$ is a random perturbation of the form
\begin{equation}\label{Hardy-condition}
	H=H^{\epsilon}:=\epsilon~H_{\epsilon}\quad\mbox{\rm with}\quad\forall n\geq 1,\quad\sup_{\epsilon \in [0,1]}
{\EE(\Vert H_{\epsilon}\Vert^n)}<\infty
\end{equation}
for some collection of centered random matrices $H_{\epsilon}$ independent of  $\Wa_t$.

\begin{theo}\label{tcl+initial-condition}
 For any $Q\in \Sa_r^+$ and any $n\geq 1$ there exists some $ \epsilon_n\in [0,1]$ such that for any $\epsilon\leq \epsilon_n$, any time horizon $t\geq 0$, and any $k=1,2$ we have 
\begin{equation}\label{estimate-tcl}
\begin{array}{l}
\vertiii{ \left[\epsilon^{-1}\left(\phi^{\epsilon}\left(Q+\epsilon~H_{\epsilon}\right)-\phi\left(Q\right)\right)\right]^k-\left[\nabla\phi(Q)\cdot H_{\epsilon}+\partial\phi\left(Q\right)\right]^k
}_{t,n}
\leq \epsilon^{1/2}~\overline{e}(t)~\vert Q\vert
\end{array}
\end{equation}

\end{theo}

The proof of this theorem is provided in Section~\ref{tcl+initial-condition-proof}. The following corollary is a direct consequence of Slutsky's lemma. 

\begin{cor}\label{cor-tcl}
Let $H_{\epsilon}$ be a sequence of processes satisfying the uniform moment condition (\ref{Hardy-condition}). In this case, we have the functional central limit theorem
$$
H_{\epsilon}\longhookrightarrow_{\epsilon\rightarrow 0}\Ha\Longrightarrow \epsilon^{-1}\left[\phi^{\epsilon}_t\left(Q+\epsilon~H_{\epsilon}\right)-\phi_t\left(Q\right)\right]\longhookrightarrow_{\epsilon\rightarrow 0} \nabla\phi_t(Q)\cdot \Ha+\partial\phi_t\left(Q\right)
$$
For any $Q\in \Sa_r^+$ and any sufficiently small $\epsilon$ we have the non-asymptotic variance estimate
\begin{equation}\label{na-var-est}
\begin{array}{l}
\Vert
\epsilon^{-2}~\EE\left(\left[\phi^{\epsilon}_t\left(Q+\epsilon~H_{\epsilon}\right)-\phi_t\left(Q\right)\right]^2\right)-
\EE\left(\left[\nabla\phi_t(Q)\cdot H_{\epsilon}+\partial\phi_t\left(Q\right)\right]^2\right)\Vert
\leq  \epsilon^{1/2}~\overline{e}(t)~\vert Q\vert
\end{array}
\end{equation}
\end{cor}

\subsection{Ensemble Kalman-Bucy filters}\label{sec-EKF}

The mathematical foundations and the convergence analysis of {\tt EnKF} algorithms are rather recent. In~\cite{legland} and~\cite{mandel}, the authors study the $\LL_n$-mean error estimates for discrete-time {\tt EnKF} algorithms. In a linear-Gaussian setting, the authors show that the {\tt EnKF} converges towards the Kalman filter as the number of samples tends to infinity. Non-linear state-space models are also considered in \cite{legland}. In~\cite{law2016,deWiljes2018} the authors consider continuous-time (non-Gaussian) state-space models (e.g. certain nonlinear diffusion models) and the convergence (with sample size) of particular {\tt EnKF} methods. 

Time-uniform fluctuation, stability and contraction estimates were given in the one-dimensional linear-Gaussian case in \cite{dbkr-2017}; with essentially no further assumptions on the underlying state-space model. In the multi-dimensional linear-Gaussian case, time-uniform stability estimates were developed under strong signal stability assumptions in \cite{delmoral16a}. Both \cite{delmoral16a,dbkr-2017} were extended in the multi-dimensional setting in \cite{Bishop/DelMoral/multiDimRicc} where assumptions on the underlying signal stability were relaxed (based on results in \cite{Bishop/DelMoral/STV2018}). Related work in \cite{tong,Kelly2014,dkt17} considers the long-time behaviour of the {\tt EnKF} in certain (possibly nonlinear) stable (and/or bounded) signal models. See also~\cite{majda-2,tong-2,Kelly2014} for a related stability analysis in the presence of (adaptive) covariance inflation and projection techniques. Typically, rather strong observability assumptions are made in all multi-dimensional {\tt EnKF} stability analysis, e.g. \cite{tong,Kelly2014,delmoral16a,dkt17,Bishop/DelMoral/multiDimRicc,deWiljes2018} and in this work.

Because of their practical importance, this subsection is dedicated to the illustration of our main perturbation analysis and fluctuation results within the {\tt EnKF} framework.

Now we introduce the relevant models in this work. Consider a time-invariant linear-Gaussian filtering (signal and observation) model of the following form
\begin{equation}\label{lin-Gaussian-diffusion-filtering}
dX_t=A\,X_t~dt+R_W^{1/2}\,dW_t\quad\mbox{\rm and}\quad
dY_t=B\,X_t~dt+R_V^{1/2}\,dV_{t}
\end{equation}
where $(W_t,V_t)$ is an $(r+r^{\prime})$-dimensional Brownian motion, $X_0$ is a $r$-valued Gaussian random vector with mean and covariance matrix $(\EE(X_0),P_0)$ (independent of $(W_t,V_t)$), the symmetric matrices $R_W$ and $R_V$ are strictly positive definite (invertible), $A$ is an arbitrary square $(r\times r)$-matrix, $B$ is an arbitrary $(r^{\prime}\times r)$-matrix, and $Y_0=0$. We let $\Fa^Y_t=\sigma\left(Y_s,~s\leq t\right)$ be the $\sigma$-algebra filtration generated by the observations. 
We now consider the conditional nonlinear McKean-Vlasov type diffusion process
\begin{equation}\label{Kalman-Bucy-filter-nonlinear-ref}
d\overline{X}_t=A~\overline{X}_t~dt~+~R^{1/2}_W~d\overline{W}_t+P_t~B^{\prime}R_V^{-1}~\left[dY_t-\left(B\overline{X}_tdt+R_V^{1/2}~d\overline{V}_{t}\right)\right]
\end{equation}
where $(\overline{W}_t,\overline{V}_t,\overline{X}_0)$ are independent copies of $(W_t,V_t,X_0)$ (independent of
 the signal and the observations). In this context, we have $S=B^{\prime}R_V^{-1}B$. Here $P_t=\Pa_{\eta_t}$ denotes the covariance matrix
\begin{equation}\label{def-nl-cov}
\Pa_{\eta_t}=\eta_t\left[(\theta-\eta_t(\theta))(\theta-\eta_t(\theta))^{\prime}\right]
\quad\mbox{\rm with}\quad \eta_t:=\mbox{\rm Law}(\overline{X}_t~|~\Fa^Y_t)\quad\mbox{\rm and}\quad
\theta(x):=x.
\end{equation}
This diffusion is a time-varying Ornstein-Uhlenbeck-type process \cite{delmoral16a,ap-2016} and thus $\eta_t$ is Gaussian.

The ensemble Kalman-Bucy filter ({\tt EnKF}) coincides with the mean-field particle interpretation of the nonlinear diffusion process \eqref{Kalman-Bucy-filter-nonlinear-ref}. To be more precise, let $(\overline{W}^i_t,\overline{V}^i_t,\overline{X}_0^i)_{1\leq i\leq N+1}$ be $(N+1)$ independent copies of $(\overline{W}_t,\overline{V}_t,\overline{X}_0)$. In this notation, the {\tt EnKF} is given by the Mckean-Vlasov type interacting diffusion process
\begin{equation}\label{fv1-3}
d\overline{X}_t^i=A~\overline{X}_t^i~dt+R^{1/2}_W~d\overline{W}_t^i+Q_tB^{\prime}R_V^{-1}\left[dY_t-\left(B ~\overline{X}_t^i~ dt+R_V^{1/2}~d\overline{V}^i_{t}\right)\right] 
\end{equation}
with $1\leq i\leq N+1$ and the rescaled particle covariance matrices  $Q_t$ defined by
\begin{equation}\label{fv1-3-2}
\begin{array}{l}
\displaystyle Q_t=Q^N_t:=\left(1+N^{-1}\right)\,\Pa_{\eta^{N}_t}
\quad\mbox{\rm with}\quad
\displaystyle\eta^{N}_t:=(N+1)^{-1}\sum_{1\leq i\leq N+1}\delta_{\overline{X}_t^i}
\end{array}\end{equation}

Observe that if $S=B^{\prime}R_V^{-1}B$ is invertible then,
\begin{equation}\label{filtering-pb-signal}
\Xa_t=S^{1/2}X_t~~\Longrightarrow ~~d\Xa_t=\overline{A}\,\Xa_t~dt+\overline{R}\,dW_t
\end{equation}
with the matrices $(\overline{A},\overline{R})$ introduced in (\ref{ref-overline-A}), and
\begin{equation}\label{filtering-pb-obs}
\Ya_t=B^{\prime}R^{-1}_VY_t\Longrightarrow d\Ya_t=\overline{B}~\Xa_t~dt+\overline{R}_V^{1/2}~dV_t\quad
\mbox{\rm with}\quad (\overline{B},\overline{R}_V)=(S^{1/2},S)
\end{equation}
In this situation we have 
$
\overline{S}:=\overline{B}^{\prime}\overline{R}^{-1}_V\overline{B}=S^{1/2}S^{-1}S^{1/2}=I
$.

From~\cite{delmoral16a}, and using the representation theorem (theorem 4.2~\cite{karatzas}; see also~\cite{doob}), there exists a filtered probability space enlargement under which we find (\ref{f21}) with the parameters
 $$
 \epsilon:={2}/{\sqrt{N}}\quad\mbox{\rm and}\quad  (R,S):=(R_W, B^{\prime}R^{-1}_VB)
 $$
That is, flow of the sample covariance (\ref{fv1-3-2}) associated with the $(N+1)$-interacting particle systems (\ref{fv1-3}), is given by the matrix Riccati diffusion (\ref{f21}) with these parameters. Thus, the stochastic Riccati equation (\ref{f21}) captures explicitly the evolution of the sample covariance of the {\tt EnKF}. The interacting particle systems (\ref{fv1-3}) are themselves interacting via the sample covariance. Note that proving $Q_t=Q^N_t\hookrightarrow_{N\rightarrow\infty} P_t$ implies roughly after some work, the convergence in some sense of $\eta^{N}_t\rightarrow\eta_t$; e.g. see \cite{delmoral16a}, and the central limit-type results given later in this section.

We note that the main results (cf. theorems \ref{theo-intro-phi-MM-Taylor}, \ref{th1-estimates-Ln}, \ref{th2} and \ref{tcl+initial-condition}) are given in terms of sufficient conditions only, and they typically require $\epsilon\leq1$ sufficiently small. In terms of the {\tt EnKF} relationship $\epsilon=2/\sqrt{N}$, this means we require $N\geq4$ sufficiently large; i.e. typically much more than $N+1=5$ particles in (\ref{fv1-3}). The sufficient conditions given in our main results may typically be stronger than required. In the next subsection, we sketch a more detailed analysis of the Riccati diffusion with $r=1$. There we suspend the requirement $\epsilon\leq1$, and we explicitly derive the limiting distribution of the diffusion. This allows us to capture tighter (necessary) conditions for moment existence, etc. 

However, we also note that in higher-dimensions (e.g. $r\geq4$), it is quite reasonable to take $N$ large enough (so that $\epsilon\leq1$ automatically) to ensure the Riccati diffusion (\ref{f21}) is well-behaved. For example, when we interpret the matrix diffusion (\ref{f21}) as the flow of some sample covariance, then intuitively one would like $N\geq r$ so the rank of the solution behaves nicely. The sufficient conditions in our main results are stronger than this anyway. In practical {\tt EnKF} applications, one typically adds some regularisation procedure \cite{evensen03,aps-2016} to ensure an associated regularised matrix diffusion is well-behaved when the number of particles is too small; but we don't discuss this idea here, see \cite{aps-2016}.

The sample mean $\eta^N_t(\theta)=:m_t=\psi_{t}^{\epsilon}(m_0,Q_0)$ is given by the stochastic flow
 \begin{equation}\label{flow-sm}
 d\psi_{t}^{\epsilon}(x,Q)=\left[A-\phi^{\epsilon}_t(Q)S\right]~\psi_{t}^{\epsilon}(x,Q)~dt+\phi^{\epsilon}_t(Q)B^{\prime}R^{-1}_V dY_t+\frac{\epsilon}{\sqrt{4+\epsilon^2}}~\Sigma_{\varphi}(\phi_t^{\epsilon}(Q))~d\overline{\Wa}_t
 \end{equation}
 with an $r$-Wiener process $\overline{\Wa}_t$ independent of  $\Wa_t$.  
Now, entering the fluctuations of the semigroup $\phi^{\epsilon}_t(Q)$ around $\phi_t(Q)$ into the flow (\ref{flow-sm})
we check that the mapping $\epsilon\mapsto\psi_{t}^{\epsilon}(x,Q)$ is differentiable with a first order derivative
$$
\partial \psi_{t}(x,Q)=\int_0^t~E_{s,t}(Q)~\left[\partial\phi_s(Q)~S^{-1/2}~\left(d\Ya_s-S^{1/2} \psi_{s}(x,Q)~ds\right)+\Sigma_{\varphi}(\phi_s(Q))~d\overline{\Wa}_s\right]
$$

Also recall that
\begin{eqnarray*}
	H^{N}_0:=\sqrt{N}~(Q_0-P_0)&\Longleftrightarrow& Q_0=P_0+\frac{1}{\sqrt{N}}~H^{N}_0\geq 0
\\
	h^{N}_0:=\sqrt{N+1}~(m_0-\EE(X_0))&\Longleftrightarrow& m_0=\EE(X_0)+\frac{1}{\sqrt{N+1}}~h^{N}_0
\end{eqnarray*}
In addition $h^{N}_0$ and $H^{N}_0$ are independent.

 We may now connect the {\tt EnKF} to our main results (cf. theorems \ref{theo-intro-phi-MM-Taylor}, \ref{th1-estimates-Ln}, \ref{th2} and \ref{tcl+initial-condition}) in a straightforward manner. For one example, by theorem~\ref{th1-estimates-Ln}, the sample path $\LL_n$-boundedness properties of the matrix diffusion process (\ref{f21}) is ensured as soon as $N>8nr$, for any $n\geq 1$. In this case, we also have $Q_t\geq 0$ for any time horizon $t\geq 0$. A time change $t\mapsto 2^{-1}\epsilon\sqrt{ N}~t~$ of the process shows that this property is also met for any real values $\epsilon\in [0, 1/\sqrt{2nr}[$. 

We also have the central limit theorem $
(h^N_0,H^N_0)\longhookrightarrow_{N\rightarrow\infty} (h_0,\Ha_0)
$,
where $h_0\stackrel{law}{=}X_0-\EE(X_0)$ and $\Ha_0$ is an independent symmetric $(r\times r)$-matrix with centered Gaussian entries equipped with a symmetric Kronecker
covariance structure
$$
\EE\left((\Ha_0\otimes \Ha_0)^{\sharp}\right) =2~(P_0\frownotimes P_0)=\EE\left((H^N_0\otimes H^N_0)^{\sharp}\right)
$$
A more detailed discussion on this multivariate central limit theorem can be found in~\cite{kendall}; see also the more recent study~\cite{wishart-bdn-17} as well as ~\cite{browne} for non-necessarily Gaussian variables. As verified later in corollary~\ref{prop-HN-sample-covariance}, the sample covariance satisfies the required moment condition (\ref{Hardy-condition}). Consequently, corollary~\ref{cor-tcl} yields the functional central limit theorem
$$
\sqrt{N}~\left[Q_t-P_t\right]~\longhookrightarrow_{N\rightarrow\infty} \nabla\phi_t(P_0)\cdot \Ha_0+\partial\phi_t\left(P_0\right)
$$
A closed form description of the Fr\'echet derivative $\nabla\phi_t(Q)$ w.r.t. the parameter $Q$
is provided in theorem~\ref{lem-frechet-derivatives}. 
In the same vein we check that
$$
\sqrt{N+1}~\left[\,m_t-\EE\left(X_t~|~\Fa_t^Y\right)\,\right]~\longhookrightarrow_{N\rightarrow\infty}\,\partial \psi_{t}(\EE(X_0),P_0)+E_t(P_0)\,h_0+\nabla\psi_t(\EE(X_0),P_0)\cdot \Ha_0
$$
In the above display $\nabla\psi_t(x,Q)\cdot H$ stands for the Fr\'echet derivative w.r.t. the parameter $Q$ given by
the formula
$$
\nabla\psi_t(x,Q)\cdot H=(\nabla E_t(Q)\cdot H)\,x+\int_0^t\,\left[
\left(\nabla E_{s,t}(Q)\cdot H\right)\,\phi_s(Q)+E_{s,t}(Q)\,\left(\nabla\phi_s(Q)\cdot H\right)\right]\,S^{-1/2}~d\Ya_s
$$
A closed form description of the Fr\'echet derivatives $\nabla E_{s,t}(Q)$  w.r.t. the parameter $Q$
is provided in corollary~\ref{cor-nabla-Pi} and its proof in Appendix \ref{appendix-section}.

In addition using (\ref{na-var-est}), for any $N$ sufficiently large we have the variance estimate
$$
\begin{array}{l}
	\Vert\,N\,\EE\left[(Q_t-P_t)^2\right]- \left(\EE\left[(\nabla\phi_t(P_0)\cdot \Ha_0)^2\right]+\EE\left[(\partial\phi_t(P_0))^2\right]\right)\Vert~\leq~  N^{-1/2}~\overline{e}(t)\,\vert P_0\vert
\end{array}
$$

We remark that a `deterministic' form of the {\tt EnKF} (abbreviated {\tt DEnKF}) introduced by Sakov and Oke in \cite{sakov2008a} yields a simpler Riccati diffusion; see also the related work \cite{Bergemann2012,Reich2013,Taghvaei2016ACC}. That Riccati diffusion will exhibit smaller random fluctuations, due to the absence of random observation perturbations in (\ref{fv1-3}); i.e. the absence of {\small$d\overline{V}^i_{t}$}, and its replacement with a deterministic adjustment. More precisely,
replace $(B\overline{X}_t+R_V^{1/2}d\overline{V}_t)$ by $B(\overline{X}_t+\eta_t(\theta))/2$ in (\ref{Kalman-Bucy-filter-nonlinear-ref}). Then, we get a special case (see \cite{Bishop/DelMoral/multiDimRicc}) of the Riccati diffusion (\ref{f21}) with $dM_t:=[Q_t^{1/2}\,d\Wa_t\,R_W^{1/2}]_{ sym}$ and $R=R_W$. This simplified Riccati diffusion captures the flow of the sample covariance for the {\tt DEnKF}. The analysis of (\ref{f21}) simplifies considerably in this situation. Naturally, we expect the Riccati diffusion in this special case to exhibit less fluctuation, and be more numerically stable (e.g. with regards to time-discretisation) because of the reduced order diffusion term. Analogous statements/simplifications can be made about the sample mean (\ref{flow-sm}). See \cite{Bishop/DelMoral/multiDimRicc} for more details. In the next subsection, we briefly study the one-dimensional Riccati diffusion with this replacement, where these properties are more readily apparent; see also \cite{dbkr-2017}. 

Going forward more broadly, we focus on the most general form of the Riccati diffusion as written in (\ref{f21}). This corresponds (ironically), to the sample covariance (\ref{fv1-3-2}) of the most naive implementation of the {\tt EnKF} given directly in (\ref{fv1-3}); i.e. without a replacement of the type just discussed. This analysis establishes a baseline for later comparison with more advanced algorithms. The main results (cf. theorems \ref{theo-intro-phi-MM-Taylor}, \ref{th1-estimates-Ln}, \ref{th2} and \ref{tcl+initial-condition}) are all given in this general setting. The broader impact of regularisation \cite{evensen03,aps-2016}, and more sophisticated {\tt EnKF} methods (e.g. \cite{evensen03,sakov2008a}), on the sample covariance flow will be discussed elsewhere; e.g. see \cite{dbkr-2017,Bishop/DelMoral/multiDimRicc}. 

As another interesting aside, we can relate the {\tt EnKF} with the Wishart process. Indeed, Wishart processes are particular instances of stochastic Riccati diffusions. They correspond to the case
$ B=0\Rightarrow
S=0$ and $\Sigma(Q)=R
$.
The term Wishart process was coined by Marie France Bru in the pioneering articles~\cite{bru1,bru2}; see also \cite{alfonsi,cuchiero,mayerhofer}.

 \subsection{Long time behavior}\label{stability-riccati-square-root}
 
The stability and the regularity properties of the Riccati equation (\ref{f21}) with $\epsilon=0$ are well understood. We return to the filtering model discussed in Section~\ref{sec-EKF}. When $R_W>0$ and $S>0$ the filtering problem associated with (\ref{filtering-pb-signal}) and (\ref{filtering-pb-obs}) is controllable and observable. These conditions ensure the existence and uniqueness of a positive-definite fixed-point $P_{\infty}$ solving the so-called algebraic Riccati equation $\Lambda(P_{\infty})=0$. In this case, the matrix difference $A-P_{\infty}S$ is asymptotically stable even when the matrix $A$ is unstable. In addition \cite{ap-2016,bd-CARE}, there exists some constant $\rho> 0$ and some function $\beta$ on $\RR_+$ such that
$\nu>0\Longrightarrow\beta(\nu)>0$, $\lim_{\nu\rightarrow 0}\beta(\nu)=0$ and for any $t\geq s\geq 0$,
\begin{equation}\label{exp-estimate-Bucy}
\Vert E_{s,t}(Q)\Vert\leq c\,\exp{\left[-\beta(\nu)~(t-s)+\rho~\nu~\Vert Q\Vert\right]}~~\Longrightarrow~~ \sup_{0\leq s\leq t}\sup_{Q\in \Sa_r^0}\Vert E_{s,t}(Q)\Vert<\infty
\end{equation}
The Riccati flow $Q\mapsto \phi_t(Q)$ is  non-decreasing positive map and for any $t\geq 0$ we have
\begin{equation}\label{monotone-prop}
 \Vert \phi_t(Q)\Vert\leq c_1~(1+\Vert Q\Vert)\quad \mbox{\rm and}\quad
\Vert \phi_t(Q)^{-1}\Vert\leq c_2~(1+\Vert Q^{-1}\Vert)
\end{equation}
The above estimates can be turned into uniform estimates w.r.t. $Q$, as soon as $t\geq \nu$.

For a more thorough discussion on the stability and the regularity properties of deterministic Riccati equations we refer to the articles~\cite{ap-2016,aps-2016,bd-CARE} and the references therein. 

The stability properties of the Riccati diffusion (\ref{f21}) are more involved. To gather some intuition of these models, we examine the one-dimensional case. The stability exposition in this subsection is explored in more depth in \cite{dbkr-2017,Bishop/DelMoral/multiDimRicc}. In this subsection only, we relax the condition $\epsilon\in[0,1]$, and simply require $\epsilon\geq0$ to be finite such that a solution to (\ref{f21}) is well defined. When $r=1$ and $R\wedge S>0$, the drift function $\Lambda=\partial F$ is the derivative of the double-well drift function
$$
	F(Q)=-\frac{S}{3}~Q~(Q-\chi^{-})~(Q- \chi^+) 
$$
with the roots
$$
\chi^{-}:=\frac{3A}{2S}-\left[
\left(\frac{3A}{2S}\right)^2+\frac{3R}{S}\right]^{1/2}< 0< \chi^+:=\frac{3A}{2S}+\left[
\left(\frac{3A}{2S}\right)^2+
\frac{3R}{S}\right]^{1/2}
$$
In this situation, (\ref{f21}) resumes to the Langevin-Riccati drift type diffusion process
\begin{equation}\label{definition-Langevin-Riccati}
dQ_t=\partial F(Q_t)~dt+
\epsilon~\sigma(Q_t)~d\Wa_t\quad \mbox{\rm with}\quad \sigma^2(Q):=Q~\left[R+SQ^2\right]
\end{equation}
Also observe that $\partial F>0$ on the open interval $]0,\chi_+[$ and
$\partial F(0)=R>0=\sigma(0)$ so that the origin is repellent and instantaneously reflecting.
In addition, the infinitesimal generator of the diffusion (\ref{definition-Langevin-Riccati}) on $]0,\infty[$ is given 
in Sturm-Liouville form by the equation
\begin{eqnarray*}
L(f)&=&\frac{\epsilon^2}{2}~\sigma^{2}~e^{V}~\partial\left(e^{-V}~\partial f\right)\quad\mbox{\rm with}\quad
V(x)=-2\epsilon^{-2}\int_{\delta}^x~\partial F(y)~\sigma^{-2}(y)~dy
\end{eqnarray*}
for any $\delta>0$. In addition, $L$ is reversible w.r.t. the probability measure $\pi$ on 
$]0,\infty[$ defined by,
\begin{equation}\label{invariant-measure-ptilde}
\pi(dx)~\propto~1_{]0,\infty[}(x)~\exp{\left[\frac{4}{\epsilon^2} \frac{A}{\sqrt{RS}}~
\tan^{-1}\left(x~\sqrt{\frac{S}{R}}\right)\right]}
\left(\frac{x}{R+Sx^2}\right)^{2\epsilon^{-2}}~\frac{1}{x(R+Sx^2)}~dx
\end{equation}
Note that in the stationary regime, one requires $\epsilon^2(n/2-1)<1$ to have existence of the $n$-th moment. Conversely, starting from some initial condition with finite $n$-th moments, this analysis implies there will be some explosion time whenever $\epsilon^2(n/2-1)\geq 1$. 

Higher-order moments even in the one-dimensional case are still troublesome. In fact, the diffusion $Q_t$ does not have any exponential moments in the stationary regime for any finite $N\geq1$ or $\epsilon>0$. In this case, for any time horizon $t\geq 0$ and for any $\nu>0$ we have
 $$
 \mbox{\rm Law}(Q_0)=\pi~~~~\Longrightarrow~~~~\EE\left(\exp{\left(\rho~\nu~\Vert \phi^{\epsilon}_{t}(Q_0)\Vert\right)}\right)=\infty
 $$
 This implies that the exponential estimate (\ref{exp-estimate-Bucy}) cannot be easily used to analyze the fluctuations of these particular stochastic Riccati diffusions.

We also remark that the heavy tailed nature of this stationary distribution implies that numerical stability may be worrisome. In the stationary regime, it is realistic to expect samples from the tails in this case, and these may be large enough and/or frequent enough to cause numerical divergence. 
 
Lastly, observe that (\ref{definition-Langevin-Riccati}) has non-globally Lipschitz coefficients. The drift is quadratic, while the diffusion has a polynomial growth of order $3/2$. It follows by \cite{hutzenthaler} that a basic Euler time-discretization may blow up, irregardless of the boundedness properties of the diffusion.

We now briefly revisit the `deterministic' version of the {\tt EnKF} in \cite{sakov2008a} (abbreviated {\tt DEnKF}). Recall the identification $(R,S):=(R_W, B^{\prime}R^{-1}_VB)$ and the modification:
$$
	\mathrm{replace}~~\left(B\overline{X}^i_t+R_V^{1/2}d\overline{V}^i_t\right)~~\mathrm{with}~~B\left(\overline{X}^i_t+\eta^N_t(\theta)\right)/2~~~\mathrm{in}~~(\ref{fv1-3}),
$$
This leads to the analysis of a simplified Riccati diffusion (\ref{f21}) with $dM_t:=[Q_t^{1/2}\,d\Wa_t\,R^{1/2}]_{ sym}$. The diffusion (\ref{f21}) with $r=1$ in this  case is reversible with respect to the following measure
\begin{equation}\label{invariant-measure-ptilde-bis}
	\pi(dx)~\propto~1_{]0,\infty[}(x)~x^{\frac{2}{\epsilon^2}-1}~\exp{\left[-\frac{S}{R\epsilon^2}~\left(x-2~\frac{A}{S}\right)^2\right]}~dx
\end{equation}
Note this measure has Gaussian tails, and we contrast this with the heavy tailed nature of (\ref{invariant-measure-ptilde}). This is significant, since it implies that the sample covariance (and mean) of this {\tt DEnKF} \cite{sakov2008a} will exhibit smaller fluctuations (than the {\tt EnKF}), and that all moments (including exponential moments) exist in this case for any choice of $N\geq1$ (or finite $\epsilon>0$). We can also expect better numerical stability (e.g. less outliers); including better time-discretisation properties \cite{hutzenthaler}. 

As previously emphasised, this work really concerns the general version of Riccati diffusion, as written in (\ref{f21}) with $dM_t:=[Q_t^{1/2}\,d\Wa_t\,\Sigma^{1/2}(Q_t)]_{ sym}$. Our main results are given in this context. In terms of the {\tt EnKF}, this means we consider the most naive implementation, as given directly by (\ref{fv1-3}) without any modification. This analysis provides a baseline for later comparison with more advanced {\tt EnKF} methods like \cite{sakov2008a} and others (e.g. with regularisation \cite{evensen03,aps-2016}). However in anticipation, we have just shown that a simple modification of the naive {\tt EnKF} (cf. the {\tt DEnKF}) may yield drastically improved sample behaviour. For example, in the naive case ({\tt EnKF}) with $r=1$, higher order moments of the sample covariance will eventually blow up for any practical $N$. Conversely, all moments exist with any $N\geq1$, after making the above \cite{sakov2008a} modification/replacement ({\tt DEnKF}). Perhaps ironically, the analysis in this latter case (cf. a simplified version of (\ref{f21})) is actually easier. See \cite{dbkr-2017,Bishop/DelMoral/multiDimRicc}.

Before we move on, consider the Wishart process again. We remark that the existence of the invariant measure (\ref{invariant-measure-ptilde}) is ensured as soon as $S\wedge R>0$ and $\epsilon<1$, without any condition on the drift $A$. When $S=0$ the existence of the stationary measure requires $A<0$. When $S=0>A=-\vert A\vert$, the stationary distribution of (\ref{f21}) is given by a Gamma distribution.
In mathematical finance, this diffusion process is called the Cox-Ingersoll-Ross process. The spectral analysis of square root diffusions is well understood. The main simplification here comes from the fact the eigenfunctions of the generator $L$ can be expressed in terms of Laguerre polynomials which form a Schauder basis of 
the Hilbert space $\LL_2(\pi)$; see e.g. exercise 330 in~\cite{dp-penev}.

\section{Some smoothness properties}\label{smoothnes-section}

\subsection{A brief review on matrix analysis and Landau-type functions}\label{review-sec}
We recall the norm equivalence formulae
\begin{equation}\label{norm-equivalence}
\Vert A\Vert_2^2:=\lambda_{ max}(A^{\prime}A)~\leq~\Vert A\Vert_F^2:= \tr(A^{\prime}A)~\leq~ r\,\Vert A\Vert_2^2
\end{equation} 
for any $(r\times r)$-matrix $A$, the above estimates are valid if we replace the spectral norm by the Frobenius norm. For the spectral norm or for the Frobenius norm, it is simple to check that
\begin{equation}\label{E-norm}
\Vert \EE(A)\Vert~\leq ~
\EE(\Vert A\Vert)\qquad\mbox{\rm and}\qquad
\vertiii{P}_{n}\leq \vertiii{\tr(P)}_{n}
\end{equation} 
where $\EE(A)$ stands for the entry-wise expectation of  some random matrix $A$ and $P$ is an $\Sa_r$-valued variable.
We recall a couple of rather well-known estimates in matrix theory. For any $(r\times r)$-square matrices $(A,B)$ by a direct application of Cauchy-Schwarz inequality we have
\begin{equation}\label{form-pq-ref}
\vert\tr(AB)\vert \,\leq\, \Vert A\Vert_F\,\Vert B\Vert_F\end{equation}
For any $P,Q\in \Sa_r^0$ we also have
\begin{equation}\label{trQ2}
\Vert P\Vert^2\leq \left(\tr\left(P\right)\right)^2\leq r~\Vert P\Vert^2
\quad\mbox{\rm
and}\quad
\lambda_{min}(P)~\tr\left(Q\right)\leq \tr\left(PQ\right)\leq \lambda_{max}(P)~\tr\left(Q\right).
\end{equation}
The r.h.s. inequality is also valid when $(Q,P)\in (\Sa_r^+\times \Sa_r)$. We check this claim
using an orthogonal diagonalization of $P$ and recalling that $Q$ remains positive semi-definite (thus with non negative 
diagonal entries). 
We also quote a local Lipschitz property of the square root function $\varphi(Q):=Q^{1/2}$ on (symmetric) definite positive matrices. For any $Q_1,Q_2\in\Sa_{r}^+$
\begin{equation}\label{square-root-key-estimate}
\Vert \varphi(Q_1)-\varphi(Q_2)\Vert \leq \left[\lambda^{1/2}_{min}(Q_1)+\lambda^{1/2}_{min}(Q_2)\right]^{-1}~\Vert Q_1- Q_2\Vert
\end{equation}
for any unitary invariant matrix norm (e.g. the spectral, or Frobenius norm); see e.g.~\cite{higham} and~\cite{hemmen}.

The transition matrix associated with a smooth flow of $(r\times r)$-matrices $A:u\mapsto A_u$ is denoted
$$
\Ea_{s,t}(A)=\exp{\left[\oint_s^t A_u~du\right]}~~\Longleftrightarrow~~ \partial_t \Ea_{s,t}(A)=A_t~\Ea_{s,t}(A)\quad\mbox{\rm and}\quad
\partial_s \Ea_{s,t}(A)=-\Ea_{s,t}(A)~A_s
$$
for any $s\leq t$, with $\Ea_{s,s}=Id$, the identity matrix. Equivalently in terms of the fundamental solution matrices
$\Ea_t(A):=\Ea_{0,t}(A)$ we have
$
\Ea_{s,t}(A)=\Ea_t(A)\Ea_s(A)^{-1}
$. We also recall that 
$$
 \left\Vert   \Ea_{s,t}(A)\right\Vert_2\,\leq \,\exp{\left(\int_s^t\rho(A_u)~du\right)}
$$
With this idea in hand, for any $s\leq t$ and $Q\in\Sa_{r}^0$ we define a key exponential semigroup of interest in this work,
\begin{equation}\label{estimate-Ea-t-bis}
E_{s,t}(Q):=\exp{\left[\oint_s^t\left(A-\phi_u(Q)\right)~du\right]}~~\Longrightarrow~~ \Vert E_{s,t}(Q)\Vert\,\leq\, \exp{\left(\rho(A)(t-s)\right)}
\end{equation}
In the further development of the article, we consider the Landau-type functions
\begin{eqnarray}
g_{\nu}(q,t)&=&c_1
\exp{\left[-\beta_1(\nu)~t+c_2~\nu~q\right]}\wedge\exp{\left[-c_3~\delta~t+c_4\right]}\quad\mbox{\rm with}\quad
\delta= \max{(-\rho(A),0)}\nonumber\\
\overline{g}_{\nu}(q,t)&=&\left[\beta_1(\nu)^{-1}\exp{\left[-\beta_2(\nu)~t+c_1~\nu~q\right]}\right]\wedge \left(c_2\exp{\left[-c_3~\delta~t\right]}~1_{\delta>0}+c_4~t^n~1_{\delta=0}\right)\nonumber\\
G_{\nu}(q,t)&=& \left[\beta_1(\nu)^{-1}\exp{\left[c_1q\nu\right]}\right]\wedge \left(c_2~1_{\delta>0}+c_3~t^n~1_{\delta=0}\right)\label{def-Landau-o}
\end{eqnarray}
with some positive constants $c_i<\infty$, some integer $n\geq 1$ and some functions $\beta_i$ such that
$\nu>0\Rightarrow\beta_i(\nu)>0$ and $\lim_{\nu\rightarrow 0}\beta_i(\nu)=0$. The constants $c_i$, the integer $n\geq 1$ and the functions $\beta_i$ may vary from line to line.
When $\nu=0$ sometimes we write 
\begin{eqnarray*}
g(t)&=&c_1~
\exp{\left[-c_2\delta t\right]}\\
 \overline{g}(t)&=&c_1\exp{\left[-c_2\delta t\right]}~1_{\delta>0}+c_3~t^n~1_{\delta=0}
\quad
\mbox{\rm and}
\quad G(t)= c_1~1_{\delta>0}+c_2~t^n~1_{\delta=0}
\end{eqnarray*}
instead of 
$
g_{0}(q,t)
$, $\overline{g}_{0}(q,t)$,  and $G_{0}(q,t)$. 
In this notation, using (\ref{exp-estimate-Bucy})  for any $\nu\geq 0$
we have the estimates
\begin{equation}\label{def-gamma}
 \Vert E_{s,t}(Q)\Vert
\leq g_{\nu} (\Vert Q\Vert,t-s),
\qquad
\Vert\Gamma_t(Q)\Vert\leq G_{\nu}(\Vert Q\Vert,t)
\quad\mbox{\rm and}\quad \Vert\Pi_t(Q)\Vert\leq \vert Q\vert_+~G_{\nu}(\Vert Q\Vert,t)
\end{equation}

\subsection{Matrix functional differentials}\label{sec-derivatives}
We let $\La(\Ta_r,\Ta_{r})$ be the set of bounded linear functional from $\Ta_r$ into itself, and equipped with the spectral or the Frobenius norm
norm. Let $\Ua\subset \Ta_r$ be an open subset. We recall that a mapping $\Upsilon:\Ua\mapsto \Ta_r$
is Fr\'echet differentiable at some $Q\in \Ua$ if there exists a 
continuous linear functional
$
\nabla \Upsilon(Q)\in  \La(\Ta_r,\Ta_{r})
$
such that
$$
\lim_{\Vert H\Vert\rightarrow 0}\Vert H\Vert^{-1}
\Vert \Upsilon(Q+H)-\Upsilon(Q)-\nabla \Upsilon(Q)\cdot H \Vert=0
$$
The mapping is said to be twice Fr\'echet differentiable at $Q\in \Ua$ when $\Upsilon$ and
the mapping
$$
\nabla \Upsilon~:~Q\in \Ua\mapsto \nabla \Upsilon(Q)\in \La(\Ta_r,\Ta_{r})
$$
is also Fr\'echet  differentiable, and so on. Given some Fr\'echet differentiable mapping $\Upsilon$ at any order at some $Q\in \Ua$, for any $H\in \Ta_r$ s.t. $Q+uH\in \Ua$ for any $u\in [0,1]$ we have
$$
\Upsilon(Q+H)=\Upsilon(Q)+
\sum_{1\leq k\leq n}~\frac{1}{k!}~\nabla^k\Upsilon(Q)\cdot H^{\otimes k}+\overline{\nabla}^{n+1}\Upsilon\left[Q,H\right]
$$
with the $(n+1)$-th order remainder
functional in the Taylor expansion given
\begin{eqnarray}\label{integral-remainder}
\overline{\nabla}^{n+1}\Upsilon\left[Q,H\right]&:=&\frac{1}{n!}~\int_0^1~(1-u)^{n}~\nabla^{n+1}
\Upsilon\left(Q+u~H\right)\cdot H^{\otimes (n+1)}~du
\end{eqnarray}
To simplify the
presentation, sometimes we write $\nabla^k\Upsilon(Q)\cdot H$ instead of $\nabla^k\Upsilon(Q)\cdot H^{\otimes k}$. 
We also consider the multilinear operator norm
$$
\vertiii{\nabla^n\Upsilon(Q)}=\sup_{\Vert H\Vert=1} 
\Vert \nabla^n\Upsilon(Q)\cdot H\Vert\quad\mbox{\rm and set}\quad
\vertiii{\overline{\nabla}^{n}\Upsilon\left[Q,H\right]}:=\Vert H\Vert^{-n}~\Vert \overline{\nabla}^{n}\Upsilon\left[Q,H\right]\Vert
$$
When $t\mapsto \Upsilon_t(Q)$  is a random process for any $m\geq 1$ and $t\geq 0$ we also set
$$
\vertiii{\nabla^n\Upsilon(Q)}_{t,m}=\sup_{\Vert H\Vert=1} 
\vertiii{ \nabla^n\Upsilon(Q)\cdot H}_{t,m}\quad\mbox{\rm and}\quad
\vertiii{\overline{\nabla}^{n}\Upsilon\left\{Q,H\right\}}_{t,m}:=\Vert H\Vert^{-n}~\vertiii{ \overline{\nabla}^{n}\Upsilon\left[Q,H\right]}_{t,m}
$$
 For any $1\leq k\leq n$ we let  $ \Ia_{k,n}$  be the set of $n$-multi-indices
\begin{eqnarray*}
 \Ia_{k,n}&:=&\{a=(a(1),\ldots,a(k))\in(\NN-\{0\})^k~:~\vert a\vert=\sum_{1\leq i\leq k}a(i)=n\}\quad
\\
\Ia_n&:=&\cup_{1\leq k\leq n}\Ia_{k,n}\quad\mbox{\rm and for any $a\in \Ia_{k,n}$ we set}\quad
 l(a)=k\quad\mbox{\rm and}\quad a!=a(1)!\ldots a(k)!.
\end{eqnarray*}
We consider the collection of matrices
\begin{eqnarray*}
\nabla^{a}\Upsilon(Q)\cdot H&:=&\left(\nabla^{a(1)}\Upsilon(Q)\cdot H,\ldots,\nabla^{a(k)}\Upsilon(Q)\cdot 
H\right)\\
\nabla^{(n,a)}\Upsilon(Q,H)&=&\left(\nabla^{(n,a(1))}\Upsilon(Q,H),\ldots,\nabla^{(n,a(k))}\Upsilon(Q,H)\right)
\end{eqnarray*}
with the matrices $\nabla^{(n,l)}\Upsilon(Q,H)$ defined for any $1\leq l\leq n$ by
$$
\displaystyle
\nabla^{(n,l)}\Upsilon(Q,H)=1_{l<n}~\frac{1}{l!}~~
\nabla^{l}\Upsilon(Q)\cdot H+1_{l=n}~~
\overline{\nabla}^{n}\Upsilon[Q,H]
$$

In this notation, we have the following Fa\`a  di Bruno's formula with remainder.
\begin{lem}\label{lem-fdb}
Consider smooth functionals $\Upsilon_1~:~\Ua\mapsto\Ta_r$, and $\Upsilon_2~:~\Va\mapsto\Ua$ on 
some open subsets $\Ua,\Va\subset \Ta_r$. For any $n\geq 1$ we have 
  \begin{equation}\label{faa-db-formula}
 \begin{array}{l}
\displaystyle ~\frac{1}{n!}~\nabla^{n}\left(\Upsilon_1\circ \Upsilon_2\right)(Q)\cdot H
\displaystyle=
\sum_{a\in \Ia_{n}}~\frac{1}{a!l(a)!}~\nabla^{l(a)}\Upsilon_1(\Upsilon_2(Q))\cdot\left(\nabla^{a}\Upsilon_2(Q)\cdot H\right)
 \end{array}
 \end{equation}
and the remainder term
  \begin{equation}\label{faa-db-formula-remainder}
\begin{array}{l}
\displaystyle\overline{\nabla}^n(\Upsilon_1\circ\Upsilon_2)[Q,H]\\
\\
\displaystyle=
\sum_{1\leq k<n}\sum_{n\leq m\leq nk}~\sum_{a\in \Ia_{k,m}}~\frac{1}{k!}~\nabla^k\Upsilon_1(\Upsilon_2(Q))\cdot \nabla^{(n,a)}\Upsilon_2(Q,H)+\overline{\nabla}^n\Upsilon_1\left[\Upsilon_2(Q),\overline{\nabla}\Upsilon_2[Q,H]\right]
\end{array}
 \end{equation}
 \end{lem}
In the same vein, we readily check the following technical lemma.
\begin{lem}
Consider smooth functionals $\Upsilon_i~:~\Ua\mapsto\Ta_r$ on 
some open subsets $\Ua\subset \Ta_r$, with $1\leq i\leq m$. The product functional
$$
\Upsilon(Q):=\Upsilon_1(Q)\ldots\Upsilon_m(Q)
$$
is also smooth on $\Ua$ with $n$-th derivatives given by the non commutative
Leibniz formula
\begin{equation}
\nabla^n \Upsilon(Q)\cdot H=\sum_{k_1+\ldots+k_m=n}~\left(\begin{array}{c}
n!\\
k_1\ldots k_m
\end{array}\right)~\nabla^{k_1} \left(\Upsilon_1(Q)\cdot H\right)\ldots \nabla^{k_m} \left(\Upsilon_m(Q)\cdot H\right)
\end{equation}
In addition, the $n$-th remainder terms are connected by the formula 
$$
\overline{\nabla}^n \Upsilon(Q)\cdot H=\sum_{n\leq k_1+\ldots+k_m\leq nm}~\nabla^{(k_1,n)} \Upsilon_1(Q,H)\ldots \nabla^{(k_m,n)}\Upsilon_m(Q,H)
$$
\end{lem}

We end this section with a Taylor expansion of square root functionals.

\begin{theo}[\cite{dp-niclas-1}]\label{lem-square-root-taylor}
The square root functional $\varphi$ is Fr\'echet differentiable at any order on $\Sa_r^+$. In addition, for any
$(A,H)\in (\Sa_r^+ \times \Sa_r)$ s.t. $A+\epsilon~H\in \Sa_r^0$ for any $\epsilon\in [0,1]$ and for any 
$n\geq 1$ 
we have
the estimates
\begin{eqnarray}
\vertiii{
\nabla^{n} \varphi(A)}\vee
\vertiii{ \overline{\nabla}^{n}\varphi\, [A,H]}&\leq& c~\lambda_{ min}(A)^{-(n-1/2)}
\label{estimate-square-root-n-nabla}
\end{eqnarray}

\end{theo}

\begin{cor}\label{cor-square-root-sigma-phi}
Let $\Upsilon$ be some smooth positive map on $\Sa^0_r$ satisfying the minorization condition
$$
\forall Q\in \Sa^0_r\qquad\Upsilon(Q)\geq R>0\quad\mbox{\rm and}\quad \vertiii{\nabla^k\Upsilon (Q)}\leq \vert
Q\vert_+ \quad \forall k\leq n
$$ 
for some $R\in \Sa_r^+$. In this situation, we have
\begin{equation}\label{estimate-varphi-Sigma}
\vertiii{\nabla^{n}\left(\varphi\circ\Upsilon\right)Q}\vee\vertiii{
\overline{\nabla}^{n}\left(\varphi\circ\Upsilon\right)[Q,H]}
\displaystyle\leq \vert Q\vert_+\max_{1\leq k\leq n }
\lambda_{ min}(R)^{-(k+1/2)}~
\end{equation}
\end{cor}
Next corollary is a direct consequence of (\ref{faa-db-formula-remainder}).
\begin{cor}\label{cor-square-root-phi}
Let $\Upsilon$ be some smooth positive map on $\Sa^0_r$ satisfying 
 for any $n\geq 1$ the following conditions
$$
\vertiii{\nabla^{n}\Upsilon(Q)} \vee \vertiii{\overline{\nabla}^{n}\Upsilon(Q,H)}\leq 
\vert H\vert_+^{u}
\vert Q\vert_+^{v}~
$$
for some $u,v\in \{0,1\}$.
In this situation, we have the estimate
\begin{equation}\label{ref-diff-varphi-phi}
\vertiii{\nabla^n(\varphi\circ\Upsilon)(Q)}\vee \vertiii{ \overline{\nabla}^n(\varphi\circ\Upsilon)(Q,H)}\leq 
\max_{1\leq k\leq n }
\lambda_{ min}(\Upsilon(Q))^{-(k+1/2)}~\vert H\vert_+^u
\vert Q\vert_+^v
\end{equation}
\end{cor}

\subsection{Semigroup expansions}

Observe that $\Lambda$
is twice Fr\'echet differentiable on $\Sa_r^0$ with first and second order 
derivatives given by
\begin{equation}\label{ref-nabla2-Lambda}
\nabla\Lambda\left(Q\right)\cdot H:=(A-Q)H+H(A-Q)^{\prime}
\quad\mbox{\rm and}
\quad
\frac{1}{2}~\partial^2 \Lambda(Q) \cdot\left[H,H\right]:=-H^2
\end{equation}
for any $H\in \Sa_r$ s.t. $Q+H\in \Sa_r^0$.
The derivatives of higher order are null. 
We also have the mean value formula
\begin{equation}
\Lambda(Q_1)-\Lambda(Q_2)=\nabla\Lambda\left(\frac{Q_1+Q_2}{2}\right)\cdot(Q_1-Q_2)\quad\mbox{\rm and}\quad
\Lambda(Q)=\nabla\Lambda(Q)\cdot Q+\Sigma(Q)
\label{def-Lambda}
\end{equation}

\begin{theo}\label{lem-frechet-derivatives}
For any $t\geq 0$  the mapping $\phi_t$ is  Fr\'echet  differentiable at any order on $\SS_r^+$.
The  derivatives are defined for any $n\geq 0$, $Q\in \SS_r^+$ and $H\in \SS_r$  by
\begin{eqnarray*}
\nabla^{n+1}\phi_t(Q)\cdot H=(-1)^{n}~(n+1)!~E_{t}(Q)~H~
(\Gamma_t(Q)~H)^n~E_{t}(Q)^{\prime}
\end{eqnarray*}
In addition whenever $Q+uH\in \SS_r^0$ for any $u\in [0,1]$,  for any $n\geq 1$ we have the estimates
\begin{eqnarray}\label{estimate-nabla-k-phi}
\vertiii{\nabla^{n}\phi_t(Q)}&\leq& \overline{g}_{\nu}\left(\Vert Q\Vert,t\right)\quad\mbox{and}\quad
\vertiii{\overline{\nabla}^{n}\phi_t\left[Q,H\right]}\leq  \overline{g}_{\nu}\left(\Vert Q\Vert+\Vert H\Vert,t\right)~
 \end{eqnarray} 

\end{theo}

\proof
By lemma 2.1 in the article~\cite{aps-2016} we have
\begin{eqnarray*}
\nabla \phi_t(Q)\cdot H&=&   E_t(Q)~H~E_t(Q)^{\prime}
\end{eqnarray*}
We consider the operator
$$
 \Delta_{t}(Q)\cdot (Q_1,Q_2):=\int_0^t~
E_{u,t}(Q)~ Q_1(u)~Q_2(u) ~E_{u,t}(Q)^{\prime}~du
$$
By proposition 3.3 in~\cite{ap-2016} we have the decomposition
$$
\phi_t(Q+H)-\phi_t(Q)~=~E_{t}(Q)~ H~E_{t}(Q)^{\prime}- \Delta_{t}(Q)~\cdot\left[\phi_t(Q+H)-\phi_t(Q),\phi_t(Q+H)-\phi_t(Q)\right]\nonumber
$$
This integral decomposition implies that
\begin{equation}\label{induction-nabla-phi}
\begin{array}{l}
\displaystyle
\nabla^n\phi_t(Q)\cdot H
\displaystyle=-\sum_{p+q=n-2}\frac{n!}{(p+1)!(q+1)!}~
\Delta_t(Q)\cdot\left[\nabla^{p+1}\phi_t(Q)\cdot H,\nabla^{q+1}\phi_t(Q)\cdot H~
\right]
\end{array}
\end{equation}
Using (\ref{exp-estimate-Bucy}), for $n=1$ we have
$$
\nabla\phi_t(Q)\cdot H=E_t(Q)~H~E_t(Q)^{\prime}~\Longrightarrow~
\vertiii{\nabla\phi_t(Q)}\leq  g_{\nu}(\Vert Q\Vert,t)
$$
We use an induction w.r.t. the parameter $n$ to prove that
$$
\begin{array}{l}
\nabla^{n+1}\phi_t(Q)\cdot H=(-1)^{n}~(n+1)!~E_{t}(Q)~H~
(\Gamma_t(Q)~H)^n~E_{t}(Q)^{\prime}
\Longrightarrow
\vertiii{\nabla^{n+1}\phi_t(Q)}\leq  \overline{g}_{\nu}(\Vert Q\Vert,t)~
\end{array}
$$
The last assertion is a consequence of (\ref{def-gamma}).
Assuming the result has been checked up to rank $n$ using the formula (\ref{induction-nabla-phi}) 
we check that
$$
~\nabla^{n+1}\phi_t(P)\cdot H
\displaystyle=(-1)^{n}~E_{t}(P)~H~
\int_0^t~\partial_u\left\{
(\Gamma_u(P)H)^{n}
\right\}~du~E_{t}(P)^{\prime}
$$
Using  (\ref{integral-remainder}) we check that
$$
\displaystyle \Vert H\Vert^{-(n+1)}~\Vert\overline{\nabla}^{n+1}\phi_t\left[Q,H\right]\Vert\leq \overline{g}_{\nu}\left(\Vert Q\Vert+\Vert H\Vert,t\right)
$$
This ends the proof of the theorem. \cqfd

We end this section with a series of corollaries on the above estimates.

\begin{cor}\label{cor-nabla-Pi}
The mappings $Q\mapsto E_{s,t}(Q)$, $\Gamma_t(Q)$ and $\Pi_t(Q)$ are smooth. In addition, for any $n\geq 1$ we have
the estimates
\begin{eqnarray}
\vertiii{\nabla^n E_{s,t}(Q)}&\leq& \overline{g}_{\nu}\left(\Vert Q\Vert,t-s\right),\qquad
\vertiii{ \overline{\nabla}^nE_{s,t}[Q,H]}\leq  \overline{g}_{\nu}\left(\Vert Q\Vert+\Vert H\Vert,t-s\right)~\label{nabla-En}\\
\vertiii{\nabla^n\Gamma_t(Q)}&\leq&G_{\nu}\left(\Vert Q\Vert,t-s\right),\qquad
\vertiii{ \overline{\nabla}^n\Gamma_t[Q,H]}\leq  G_{\nu}\left(\Vert Q\Vert+\Vert H\Vert,t-s\right)\label{nabla-Gamma}
\end{eqnarray}
as well as
\begin{equation}
 \vertiii{\nabla^n\Pi_t(Q)}\leq G_{\nu}\left(\Vert Q\Vert,t\right)~\vert Q\vert_+~~\mbox{and}~~~
  \vertiii{\overline{\nabla}^{n}\Pi_t\left[Q,H\right]}\leq G_{\nu}\left(\Vert Q\Vert+\Vert H\Vert,t\right)~\vert H\vert_+
\vert Q\vert_+\label{estimate-nabla-Pi}
\end{equation}
\end{cor}

The proof of the above corollary is provided in Appendix \ref{appendix-section}.

\begin{cor}
We further assume that $H$ is a random perturbation of the form
\begin{equation}\label{Hardy-condition-2}
H=H^{\epsilon}:=\epsilon~H_{\epsilon}\quad\mbox{with}\quad \sup_{\epsilon\in [0,1]}{
\EE\left(\exp{\left[\delta\Vert H_{\epsilon}\Vert\right]}
\right)}<\infty
\quad\mbox{\rm
for some $\delta>0$.}
\end{equation}
In this situation there exists some $\nu\geq 0$ such that 
$$
\phi_t(Q+\epsilon~H_{\epsilon})=\phi_t(Q)+\sum_{1\leq k< n}~(-1)^{k-1}~\epsilon^{k}
\left[E_{t}(Q)H_{\epsilon}\left(\Gamma_t(Q)H_{\epsilon}\right)^{k-1}E_{t}(Q)^{\prime}\right]+\overline{\nabla}^{n}\phi_t(Q, \epsilon~H_{\epsilon})
$$
with a stochastic remainder term such that
$$
\Vert\overline{\nabla}^{n}\phi_t(Q,\epsilon~H_{\epsilon})\Vert_{t,m}\leq \epsilon^n~\overline{g}_{\nu}\left(
\Vert Q\Vert,t\right)
$$
\end{cor}

The prototypical perturbation model we have in mind is the sample covariance matrix defined in (\ref{fv1-3-2}). We end this section around this theme. To clarify the presentation, we set 
$(X,Q,H_N)$ instead of $(X_0,P_0,H^{N}_0)$.
We let $\Pa_n(m)$ be the set of all partitions of $[n]=\{1,\ldots,n\}$ with no more than $m$ blocks. 
 We also consider the mapping
$$
\alpha^{\pi}=\sum_{1\leq i\leq \vert \pi\vert }i~1_{\pi_i}
$$
where $\pi_1,\ldots,\pi_{\vert \pi\vert}$ stand for the $\vert \pi\vert$ blocks of a partition
of $[n]$ ordered w.r.t. their smallest element. We also set
$$
\Gamma_{\pi}(Q):=~\EE\left(\left(\XX_{\alpha^{\pi}(1)}-Q\right)\prod_{2\leq l\leq \{\pi\}}\left(\Gamma_t(Q)\left(\XX_{\alpha^{\pi}(l)}-Q\right)\right)\right)
$$
with independent copies $\XX_i$ 
of the random matrix $\XX=(X-\EE(X))(X-\EE(X))^{\prime}$, and 
$$
\Pi_N^n:=\cup_{2\leq k<n}\Pa_k\left([k\wedge N]/2\right)\quad\mbox{\rm and}\quad
\{\pi\}:=k\Longleftrightarrow\pi\in\Pa_k\left((k\wedge (N-1))/2\right)
$$
Following the matrix moment expansions developed in~\cite{wishart-bdn-17} we check the following corollary. 
\begin{cor}\label{prop-HN-sample-covariance}
For any  $n\geq 1$  we have the uniform estimates
$$
\sup_{N\geq 1}{
\EE\left(\Vert H_{N}\Vert_F^n
\right)}^{1/n}\leq c~n~\tr(Q)
\Longrightarrow
\Vert \overline{\nabla}^n\phi\left[Q,N^{-1/2}~H_N\right]\Vert_{t,m}\leq \frac{1}{N^{n/2}}~\overline{g}_{\nu}\left(
\Vert Q\Vert,t\right)
$$
for some $\nu\geq 0$.
In addition, we have the non asymptotic Taylor expansion
\begin{equation}\label{non-asymptotic-Taylor-phi-QN}
N^{n/2}~\Vert\EE\left(\phi_t(Q^N)\right)-\phi_t(Q)+
\sum_{\pi\in \Pi_N^n}~(-2)^{\{\pi\}}~N^{-\{\pi\}}~(N)_{\vert\pi\vert}~\Gamma_{\pi}(Q)\Vert
\leq \overline{g}_{\nu}\left(
\Vert Q\Vert,t\right)
\end{equation}
  \end{cor}

\section{Stochastic matrix integration} \label{stoch-matrix-calc-section}

\subsection{Some algebraic aspects of tensor products}
For any $n\geq 1$ we set $[n]:=\{1,\ldots,n\}$ and $\NN_n= (\NN-\{0\})^n$.
More generally, for any multi-indices $p=(p_1,\ldots,p_n)\in \NN_n$, for some $n\geq 1$,
we set $[p]=[p_1]\times\ldots\times[p_n]$.   With a  slight abuse notation, when there are no confusion we write
$1$ instead $(1,\ldots,1)$ the $n$-multi-index with unit entries. For any $n$-multi-index $p=(p_1,\ldots,p_n)\in\NN_n$, for some $n\geq 1$, we write
 $1\leq i\leq p$ for  the set of 
multi-indices $i=(i_1,\ldots,i_n)$ s.t. $1\leq i_l\leq p_l$, with $1\leq l\leq n$.

We let $\Ta_{p,q}$ denote the tensor space spanned by all 
linear transformations from the tensor space
$\RR^{q}:=\RR^{q_1}\otimes\ldots \otimes\RR^{q_n}$
into $\RR^{p}:=\RR^{p_1}\otimes\ldots \otimes\RR^{p_m}$, for some $m,n\geq 1$.  With a slight abuse of notation we write $0$ and $I$ the null and the unit tensors in $\Ta_{p,q}$, for any multi-indices $p,q$.

In this notation, the transposition $T^{\prime}\in \Ta_{q,p} $
of a tensor $T\in \Ta_{p,q}$, and
the trace of a tensor
$\overline{T}\in \Ta_{p,p}$ are defined for any $1\leq i\leq q$ and $1\leq j\leq p$  by the same formulas
as for conventional matrices; that is we have that
$$
T^{\prime}_{i,j}=T_{j,i}\quad\mbox{\rm and}\quad
\tr(\overline{T})=\sum_{1\leq i\leq p}\overline{T}_{i,i}
$$
It is also convenient to consider the $q$-partial trace contraction $T^{\flat(p,r)}\in \Ta_{p,r}$ of a tensor
$T\in \Ta_{(p,q),(q,r)}$ defined for any $1\leq k\leq p$ and $1\leq l\leq r$ by
$$
T^{\,\flat(p,r)}_{k,l}=\sum_{1\leq i\leq q}T_{(k,i),(i,l)}
$$
When $p=r=1$ we write $T^{\, \flat}$ instead of $T^{\, \flat(1,1)}$.
We denote by $\Sa_{p,p}$ the space of symmetric tensors $\overline{T}=\overline{T}^{\,\prime}$, 
$\Sa_{p,p}^0$ and $\Sa_{p,p}^+$ the subspace of positive semidefinite and positive definite tensors.
Note the $\bullet$-tensor
product $T\,\overline{T}$ of a tensor $T\in \Ta_{p,q}$
with a tensor $\overline{T}\in \Ta_{q;r}$ is a $\Ta_{p;r}$-tensor with $(i,j)$-entries given by the
formula
$$
(~T\,\overline{T}~)_{i,j}=\sum_{1\leq k\leq p}~T_{i,k}~\overline{T}_{k,j}
$$
The  
products $T\otimes \overline{T}$ and $T~\overline{\otimes}~ \overline{T}$ of tensors $(T,\overline{T})\in (\Ta_{p,q}\times \Ta_{r,s})$ are defined by
$$
(T\otimes \overline{T})_{(i,j),(k,l)}=T_{i,k}~\overline{T}_{j,l}=
(T~\overline{\otimes}~ \overline{T})_{(i,j),(l,k)}
$$
We also write $T^{\otimes n}=T\otimes\ldots\otimes T$ the $n$-fold tensor product of a given tensor $T$. 
We let $(p,q)^{\,\sharp}$ be the transposition of a 
$(2n)$-multi-index  $(p,q)=(p_1,\ldots,p_{n},q_1,\ldots,q_n)$ by
$$
(p,q)^{\,\sharp}=((p_1,q_{1}),(p_2,q_{2})\ldots,(p_n,q_{n}))\in \NN_{2n}
$$
For $n=2$ observe that $ ((p,q)^{\,\sharp})^{\,\sharp}=(p,q)$. 
The $\sharp$-transpose $T^{\,\sharp}$ of a tensor $T\in \Ta_{(p,q)}$ is the tensor
with entries 
$ T^{\,\sharp}_{(i,j)^{\,\sharp}}=T_{(i,j)}
$.

We illustrate the tensor product properties discussed above with matrices and vectors. 

 In the further
development $A,B,\ldots$ and $A_n$ stands for any matrices with appropriate dimensions so that the formulae
make sense. We also identify $(p_1\times 1)$-matrix with $p_1$-column vectors, and more generally 
$(p,1)$-tensors with $p$-tensors for any multi-indices $p$ and $1=(1,\ldots,1)$ for any fold product.
For instance, when $q_1=1$ sometimes we write $A_{i_1}$ instead of $A_{i_1,1}$ the $i_1$-th entry of the
$p_1$-column vector $A$. 
For multi-indices $(q_1,q_2)=(1,1)$ 
we also write  $A_{(i_1,i_2)}$ instead of $A_{(i_1,i_2),1}$ 
the $(i_1,i_2)$-entry of the $(p_1,p_2)$-tensor $A$. We use the letters $x,y$ to denote column vectors.
The $\otimes$-tensor product $(A\otimes B)$ and $(A\otimes B)^{\,\sharp}$ are defined by
\begin{equation}\label{def-sharp}
(A\otimes B)_{(i_1,i_2),(j_1,j_2)}:=A_{i_1,j_1}~B_{i_2,j_2}:=(A\otimes B)^{\,\sharp}_{(i_1,j_1),(i_2,j_2)}
\quad
\mbox{\rm and}\quad
\left(A_1\otimes\ldots\otimes A_n\right)^{\flat}=A_1\ldots A_n
\end{equation}
We also have  the $\sharp$-transposition rules
\begin{equation}\label{sharp-transposition-tensor-product}
\left[(A\otimes \overline{A})~T~(B\otimes \overline{B})\right]^{\,\sharp}=
(A\otimes B^{\prime})~T^{\,\sharp}~(\overline{A}^{\prime}\otimes \overline{B})
\end{equation}

The symmetric Kronecker and the commutative Lyapunov products of $A$ and $B$ 	are defined 
 by
\begin{equation}\label{def-sym-Kr}
4(A\frownotimes B)
\displaystyle:=
(A\otimes B)+(B\otimes A)+(A\overline{\otimes}B)+(B\overline{\otimes}A)\quad\mbox{\rm and}\quad
2(A\odot B):=\displaystyle A\{B\}+\{B\}A
\end{equation}
with the matrix bracket $A\mapsto\left\{A\right\}$ defined in (\ref{def-A-par}). Observe that
\begin{equation}\label{def-sym-Kr-bis}
\begin{array}{l}
(A\otimes B)^{\sharp}_{(i,j),(k,l)}=(A\otimes B)_{(i,k),(j,l)}
\quad \mbox{\rm and}\quad
(A\frownotimes B)^{\sharp}_{(i,j),(k,l)}=(A\frownotimes B)_{(i,k),(j,l)}\\
\\
\Longrightarrow
\left((A\otimes B)^{\sharp}\right)^{\flat}=(A\otimes B)^{\flat}= AB\quad \mbox{\rm and}\quad
\left((A\frownotimes B)^{\sharp}\right)^{\flat}=(A\frownotimes B)^{\flat}=A\odot B
\end{array}
\end{equation}
We also have
\begin{equation}\label{AstarB}
P,Q\geq 0\Longrightarrow 0\leq P\odot Q\leq 2^{-1}\left(\Vert P\Vert_2~\tr(Q)+\Vert Q\Vert_2~\tr(P)\right)~I
\end{equation}
The symmetric Kronecker product is commutative but not associative. the above products are distributive w.r.t. the addition of matrices.
We recall that the Kronecker product is associative, distributive over ordinary matrix addition, and compatible
w.r.t. the ordinary multiplication of matrices. We also have the commutation properties
\begin{equation}
(A\otimes B)(C~\overline{\otimes}~D)=(AC)~\overline{\otimes}~(BD)~=~(A~\overline{\otimes}~B)~(D\otimes C)
\label{prod-prop}
\end{equation}
as well as the norm estimates
\begin{equation}\label{norm-estimates-ref-frown}
\begin{array}{l}
\Vert (A\frownotimes B)^{\,\sharp} 
\Vert\vee\Vert A\frownotimes B\Vert\leq  \Vert A\otimes B \Vert= \Vert A \Vert~\Vert B \Vert\quad\mbox{\rm and}\quad \Vert A\overline{\otimes} B\Vert_F=
 \Vert A\Vert_F\Vert B\Vert_F\\
\\
\Longrightarrow
\Vert (A\frownotimes B)^{\,\sharp}-(\overline{A}\frownotimes \overline{B})^{\,\sharp}\Vert\leq c~
\left[~\Vert B\Vert~\Vert A-\overline{A}\Vert+\Vert \overline{A}\Vert~\Vert B-\overline{B}\Vert~\right]
\end{array}
\end{equation}

For further discussion on the Kronecker product, we refer to Graham~\cite{graham} and Van Loan~\cite{loan}.

\subsection{Matrix valued martingales}

The martingales discussed in this section are to be understood w.r.t. some continuous filtration $\Fa=(\Fa_t)_{t\geq 0}$ of $\sigma$-fields on some probability space $(C([0,\infty[,\RR^{d}),\PP)$, where $C([0,\infty[,\RR^{d})$  stands for the 
set of  continuous trajectories on $\RR^d$, for some dimension parameter $d\geq 1$. A 
multivariate martingale is a matrix valued and continuous 
stochastic process whose entries are martingales. 

We recall that the angle bracket 
between a couple of real valued martingales $(A,B)$ is the differentiable  
and increasing process $$t\in [0,\infty[\mapsto
\langle A\,|\, B\rangle_t\in [0,\infty[\quad
\mbox{\rm such that}\quad
A_tB_t-\langle A\,|\, B\rangle_t\quad\mbox{\rm is a martingale}
$$
In the further development of this section we  implicitly assume that the processes are chosen so that the angle brackets of the martingales are integrable (or equivalently the martingales are square integrable).
We use the notation $\langle \point\,|\, \point\rangle$ to differentiate the angle bracket with the inner product
in $\RR^n$ and the Frobenius inner product.
For any multi-indices $p=(p_1,\ldots,p_n)$, $q=(q_1,\ldots,q_n)$, $r=(r,\ldots,r_n)$, and
 $s=(s_1,\ldots,s_n)$
   the $\otimes$-angle bracket  between an 
 a
$\Ta_{p;r}$-valued martingale $A$  and
a $\Ta_{q,s}$-valued martingale $B$ is the $\Ta_{(p,q),(r,s)}$-tensor $\langle A\,\vert\otimes \,\vert\, B\rangle$ such that
 $$A\otimes B-\langle A\,\vert\otimes \,\vert\, B\rangle\qquad\mbox{\rm is a $\Ta_{(p,q),(r,s)}$-valued martingale.}$$

When $n=1$, for any $(p\times r)$-matrix valued martingale $A$ and any $(q\times s)$-matrix valued martingale $B$, recalling that the trace of a martingale is a martingale, we have commutation formula
\begin{equation}\label{commutation-mb}
\tr(A\otimes B)=\tr(A)\,\tr(B)
\quad\Longrightarrow\quad
\langle \tr(A)\,|\, \tr(B)\rangle=
\tr\left(\langle A\,\vert\otimes \,\vert\, B\rangle\right)
\end{equation}

In multidimensional settings, another important notion is the $\bullet$-angle bracket $\langle A\,|\,B\rangle$  between 
and $(p\times n)$-matrix valued 
continuous martingale $A$ and and $(n\times q)$-matrix valued 
continuous martingale $B$ 
defined as
the $(p\times q)$ matrix $\langle A\,|\,B\rangle$ such that
$$
AB-\langle A\,|\,B\rangle \quad \mbox{\rm is a martingale}\Longrightarrow \langle A\,\vert\,\otimes\,\vert\,B\rangle^{\flat}=\langle A\,|\,B\rangle=
\left( \langle A\,\vert\,\otimes\,\vert\,B\rangle^{\sharp}\right)^{\flat}
$$

\subsection{Stochastic matrix integration}\label{section-smi}
In the further
development $A,\overline{A},B,\overline{B},M,\overline{M},P,\overline{P}, Q,\overline{Q}\ldots$  stands for any matrix valued
processes with appropriate dimensions so that the formulae
make sense.
We also use the letters $P,Q,\overline{P},\overline{Q}$ to denote positive-semidefinite matrix-valued
stochastic processes, $M,\overline{M}$ to denote martingales, and $A,\overline{A},B,\overline{B}$ any stochastic process. We also consider the stochastic matrix integrals
  \begin{equation}\label{def-AMB-bullet}
(  A\bullet M\bullet B)_t:=\int_0^t A_s\,dM_s\,B_s \end{equation}
The angle brackets between martingales of the form (\ref{def-AMB-bullet}) are given by the formula
       \begin{equation}\label{prop-ra-lg-ll-gg-sharp}
    \langle A\bullet M\bullet B^{\prime}\,\vert\,\otimes
\,\vert\, \overline{A}^{\prime}\bullet \overline{M}\bullet \overline{B}\rangle^{\,\sharp}=
(A\otimes B)\bullet\,\langle M\,\vert\, \otimes\,\vert\,\overline{M}\rangle^{\,\sharp}\bullet(\overline{A}\otimes \overline{B})
        \end{equation}
For instance, using the commutation property (\ref{prod-prop}) we check that
\begin{equation}\label{coupligng-eq-W}
\begin{array}{l}
\partial_t\,\langle \Wa\,\vert\, \otimes\,\vert\,\Wa\rangle^{\,\sharp}=I\otimes I\quad\mbox{and}\quad
 \partial_t\,\langle \Wa\,\vert\, \otimes\,\vert\,\Wa^{\prime}\rangle^{\,\sharp}=I~\overline{\otimes}~ I\\
 \\
 \Longrightarrow
 \partial_t\,
     \langle (A\bullet \Wa\bullet B^{\prime})_{ sym}\,\vert\,\otimes
\,\vert\, (\overline{A}\bullet \Wa\bullet \overline{B}^{\prime})_{ sym}\rangle^{\,\sharp}=
(A\overline{A}^{\prime})\frownotimes (B \overline{B}^{\prime})   
 \end{array}
\end{equation}
By Doob's representation theorem (see theorem 4.2~\cite{karatzas}, and the original work by Doob~\cite{doob}), the above formula shows that any symmetric continuous martingale
has an angle bracket with a symmetric Kronecker structure. More precisely, there exists some filtered probability space underwhich
$$
\begin{array}{l}
M:=(A\bullet \Wa\bullet B^{\prime})_{ sym}\stackrel{ law}{=}(Q_A^{1/2}\bullet \Wa\bullet Q_B^{1/2})_{ sym}\Longrightarrow  \partial_t\,\langle M\,\vert\,\otimes\,\vert\,M\rangle^{\,\sharp}=Q_A\frownotimes Q_B
 \end{array}$$
with the positive map $A\mapsto Q_A:=AA^{\prime}$.
We further assume that $M$ is an an $\Sa_r$-valued martingale such that
\begin{equation}\label{tr-k-mm}
\begin{array}{l}
 \partial_t\,\langle M\,\vert\,\otimes\,\vert\,M\rangle^{\,\sharp}= 
Q\frownotimes \overline{Q}
\Longrightarrow \partial_t\,\langle M\,\vert\, M\rangle=Q\odot \overline{Q}\quad\mbox{and}\quad
\partial_t\,\langle \tr(M)\,\vert\,\tr(M)\rangle=\tr\left(Q\overline{Q}\right)
\end{array}
\end{equation}
In this situation we have the angle bracket formula
\begin{equation}\label{tr-k-mm-otimes-2}
\begin{array}{l}
\displaystyle 4~ \partial_t\langle A\bullet M\bullet B^{\prime}\,\vert\, \otimes\,\vert\, \overline{A}^{\prime}\bullet M\bullet \overline{B}\rangle_t^{\sharp}\\
\\
\displaystyle=
(AQ\overline{A})\otimes
(B\overline{Q}\,\overline{B})+(A\overline{Q}\,\overline{A})\otimes(
BQ\overline{B})+
(AQ\overline{B})~\overline{\otimes}~(B\,\overline{Q}\,\overline{A})+
(A\overline{Q}\,\overline{B})~\overline{\otimes}~(B\,Q\,\overline{A})\\
\\
\Longrightarrow   \partial_t\langle A\bullet M\bullet A^{\prime}\,\vert\, \otimes\,\vert\, A\bullet M\bullet A^{\prime}\rangle_t^{\sharp}=(AQA^{\prime})\frownotimes(A\overline{Q}A^{\prime})
\end{array}
\end{equation}
  We also have the contraction angle bracket formula
\begin{equation}\label{tr-k-mm-3}
\begin{array}{l}
\displaystyle 4\partial_t\langle A\bullet M\bullet B^{\prime}\,\vert\,  \overline{A}^{\prime}\bullet M\bullet \overline{B}\rangle\\
\\
\displaystyle=\left(AQ\overline{A}\right)\left(B\,\overline{Q}\,\overline{B}\right)+
\left(A\overline{Q}\,\overline{A}\right)\left(B\,Q\,\overline{B}\right)
 \displaystyle+\left(AQ\overline{B}\right)~
 \tr\left(\overline{Q}\,\overline{A}B \right)
 \displaystyle+\left(A\overline{Q}\,\overline{B}\right)~ \tr\left(Q\overline{A}B \right)
\end{array}
\end{equation}
and well as the trace angle bracket 
\begin{equation}\label{tr-k-mm-trace-2}
\displaystyle\partial_t\,\left\langle \tr\left[  A\bullet M\bullet B \right]\,\vert\,  \tr\left[
\overline{A}\bullet M\bullet\overline{B}\right]\right\rangle=2^{-1}~\tr\left((\overline{B}\,\overline{A})_{ sym}~\left[\overline{Q}\,(BA)_{sym}\,Q+
Q\,(BA)_{sym}\,\overline{Q}\right]\right)
\end{equation}
Notice that
\begin{eqnarray*}
\displaystyle\partial_t\langle A\bullet M\bullet B^{\prime}\,\vert\,  B\bullet M\bullet A^{\prime}\rangle&=&AB^{-1}~\left(Q_B\odot\overline{Q}_B\right)~\left(AB^{-1}\right)^{\prime}\\
\displaystyle\partial_t\,\left\langle \tr\left[  A\bullet M\bullet B \right]\,\vert\,  \tr\left[
A\bullet M\bullet B\right]\right\rangle&=&\tr\left(\overline{Q}\,(BA)_{sym}\,Q\,(BA)_{ sym}~\right)
\end{eqnarray*}
with the matrices 
$
Q_B:=BQB^{\prime}$ and $ \overline{Q}_B:=B\,\overline{Q} B^{\prime}$.
We end this section with a functional Ito formula associated with the stochastic Riccati equation (\ref{f21}).

Let $\Upsilon$ be a smooth mapping from $\Sa_r$ into itself.
We let $e_i$ be the $r$-row unit vector with entries 
$e_i(k)=1_{i}(k)$, for any $1\leq k\leq r$. The set $\Sa_r$ is spanned by the symmetric matrices
$e_{i,j}=(e_i\otimes e_j)_{ sym}=(e^{\prime}_ie_j+e^{\prime}_je_i)/2$. Thus we can identify $\nabla \Upsilon(Q)$ with the tensor
$$
\nabla \Upsilon(Q)_{(i,j),(k,l)}:=
\left(\nabla \Upsilon(Q)\cdot e_{i,j}\right)_{k,l}=e_k~\left(\nabla \Upsilon(Q)\cdot e_{i,j}\right)~e_l^{\prime}
$$
In this notation, we have
$$
\begin{array}{l}
H=\sum_{i,j}~H_{i,j}~e_{i,j}
\Longrightarrow\left(\nabla \Upsilon(Q)\cdot H\right)_{(k,l)}=\sum_{i,j}~H_{(i,j)}~\nabla \Upsilon(Q)_{(i,j),(k,l)}=\left(
H~\nabla \Upsilon(Q)\right)_{(k,l)}
\end{array}$$
In the same vein,  we can identify $\nabla^2 \Upsilon(Q)$ with the tensor
$$
\begin{array}{l}
\nabla^2 \Upsilon(Q)_{((i,j),(k,l)),(m,n)}:=
\left(\nabla^2 \Upsilon(Q)\cdot (e_{i,j},e_{k,l})\right)_{m,n}
\\
\\
\Longrightarrow\left(\nabla^2 \Upsilon(Q)\cdot (H_1,H_2)\right)=(H_1\otimes H_2)^{\sharp}~\nabla^2 \Upsilon(Q)
\end{array}$$
We  let $L$ the second order differential functional
$$
L_{\epsilon}(\Upsilon)(Q)=\Lambda(Q)\nabla\Upsilon(Q)+\frac{\epsilon^2}{2}~(Q_t\frownotimes \Sigma(Q_t))\nabla^2\Upsilon(Q)
$$
In this notation we readily check the following proposition.
\begin{prop} 
For any smooth mapping $\Upsilon$ from $\Sa_r$ into itself with polynomial growth 
derivatives, and for any sufficiently small $\epsilon$ we have the Ito formula
\begin{equation}\label{ito-formula}
d\Upsilon(Q_t)=L_{\epsilon}(\Upsilon)(Q_t)~dt+\epsilon ~dM_t(\Upsilon)
\end{equation}
with the martingale $dM_t(\Upsilon)=dM_t~\nabla\Upsilon(Q_t)$ with angle bracket
$$
\partial_t\langle M(\Upsilon)\,\vert\,\otimes\,\vert\,M(\Upsilon)\rangle_t^{\sharp}=(Q_t\frownotimes \Sigma(Q_t))\left(\nabla\Upsilon(Q_t)\otimes \nabla\Upsilon(Q_t)\right)
$$
\end{prop}
\subsection{A martingale continuity theorem}

\begin{prop}\label{propM-bounds-ref}
For any $n\geq 1$ and any flows $P_s^{(i)}$ and $ Q_s^{(i)}$, $1\leq i\leq m$, we have
\begin{equation}\label{ref-M-bound-mn}
\displaystyle M:=\sum_{1\leq i\leq m}(P^{(i)}\,\bullet \,\Wa\,\bullet Q^{(i)})_{ sym}
\displaystyle\Longrightarrow\vertiii{M}_{t,n}\leq c~\sqrt{t}~ \vertiii{P}_{t,2n}~\vertiii{Q}_{t,2n}
\end{equation}
In addition, for any flows $A_s,B_s$ we have
\begin{eqnarray}
\displaystyle\vertiii{ A\bullet M\bullet B}_{t,2n}^2
\displaystyle&\leq& c~\sum_{1\leq i,j\leq m}\EE\left[\left(\int_0^t~\Vert A_s\Vert^2~\Vert B_s\Vert^2~\Vert P_s^{(i)}\Vert^{2}~
\Vert Q_s^{(j)}\Vert^{2}
~ds\right)^n\right]^{1/n}\label{ref-M-bound-mn-2-integral}\\
&\leq& c~\sqrt{t}~\vertiii{A}_{t,8n}\vertiii{B}_{t,8n}~ \vertiii{P}_{t,8n}~\vertiii{Q}_{t,8n}\label{ref-M-bound-mn-2}
\end{eqnarray}

\end{prop}
The proof of this proposition is technical, so it appears in Appendix \ref{appendix-section}. We use the same notational conventions as in Section~\ref{section-smi}. We also set
     \begin{eqnarray}
    \chi(P,\overline{P},Q,\overline{Q})&:=&c~\lambda_{min}(\overline{P},\overline{Q})^{-1} \left[1+\tr(\overline{P})^2+ \tr(P)^2+
\tr(
Q)^2+\tr(
\overline{Q})^2
\right]\nonumber\\
  \lambda_{min}(\overline{P},\overline{Q})&:=&
  1\wedge\lambda_{min}(\overline{P})\wedge \lambda_{min}(\overline{Q})\quad\mbox{\rm and}\quad
  \chi_n(t):=\vertiii{ \chi(P_t,\overline{P}_t,Q_t,\overline{Q}_t)}_{n}
  \label{def-chi-PQ}
   \end{eqnarray}
  In the further development of this section, we let $(M,\overline{M})$ be a couple of martingales of the form
\begin{equation}\label{def-MM-cor}
M:=\left(P^{1/2}\bullet \Wa\bullet Q^{1/2}\right)_{ sym}\quad\mbox{\rm and}\quad
 \overline{M}:=\left(\overline{P}^{1/2}\bullet \Wa\bullet \overline{Q}^{1/2}\right)_{ sym}
\end{equation}
   The main objective of this section is to prove the following theorem.
\begin{theo}\label{theo-bdg}
For  any time horizon $t$  any $n\geq 2$ and any stochastic processes $(A,B)$
we have the functional estimate
\begin{equation}\label{pr-fcv-2}
\begin{array}{l}
\displaystyle\vertiii{
A\bullet (M-\overline{M})\bullet B}_{t,n}^{2}
\displaystyle\leq 
\int_0^t~\chi_{2n}(s)~
\vertiii{A_s\otimes B_s}_{2n}^2~\left[~\vertiii{Q_s-\overline{Q}_s}_{2n}+\vertiii{
   P_s-\overline{P}_s}_{2n}~\right]
 ~ds
 \end{array}
\end{equation}
\end{theo}
\begin{cor}
Let $(M,\overline{M})$ be the couple of martingales defined in (\ref{def-MM-cor}).  For any $n\geq 2$ we have
\begin{equation}\label{pr-fcv}
\vertiii{ M-\overline{M}}_{t,n}^2\leq 
\int_0^t~\chi_n(s)~\left[\vertiii{
Q_s-\overline{Q}_s}_n+\vertiii{
  P_s-\overline{P}_s}_n\right]~ds
\end{equation}
\end{cor}

The proof of the theorem is based on a couple of technical lemmas of their own interest.

\begin{lem}\label{lem-pre-cv}
Under the assumptions of theorem~\ref{theo-bdg}, we have the Lipschitz estimate
 \begin{equation}\label{local-Lip}
  \displaystyle \Vert\partial_t\,\langle\,(M-\overline{M})_{ sym}\,\vert\,\otimes
\,\vert\, (M-\overline{M})_{ sym}\rangle^{\,\sharp}_t\Vert
   \displaystyle\leq     \chi(P_t,\overline{P}_t,Q_t,\overline{Q}_t)~
 \left[\Vert Q_t-\overline{Q}_t \Vert+\Vert P_t-\overline{P}_t \Vert\right]
    \end{equation}
\end{lem}

The proof of the above technical lemma is provided in Appendix \ref{appendix-section}. Next lemma is a direct consequence of
Burkholder-Davis-Gundy inequality.

\begin{lem}[Burkholder-Davis-Gundy inequality]\label{lem-bdg-inequality}
For any $t\geq 0$, and $n\geq 1$ we have
\begin{equation}\label{bdg-inequality}
  \begin{array}{l}
       \displaystyle
\vertiii{A\bullet M\bullet B}_{t,2n}^{2n}
   \displaystyle\leq 
\EE\left[\left(\int_0^t\tr\left[\left((A_s^{\prime}A_s)\otimes (B_sB_s^{\prime})\right)~
\partial_s\langle M\,\vert\,\otimes\,\vert \,M\,\rangle^{\,\sharp}_s
\right]~ds\right)^n\right]
    \end{array}
\end{equation}

\end{lem}
The proof of theorem~\ref{theo-bdg} is now almost immediate.

{\bf Proof of theorem~\ref{theo-bdg}:} Combining (\ref{bdg-inequality}) with lemma~\ref{lem-pre-cv} and using the generalized Minkowski inequality for any $n\geq 1$ we have
$$
 \begin{array}{l}
     \displaystyle
\vertiii{M_s-\overline{M}}_{t,2n}^2
   \displaystyle\leq 
\int_0^t~\EE\left[\chi(P_s,\overline{P}_s,Q_s,\overline{Q}_s)^n
 \left[\Vert Q_s-\overline{Q}_s \Vert_F+\Vert P_s-\overline{P}_s \Vert_F\right]^n\right]^{1/n}~ds
    \end{array}
$$
The end of the proof of the first assertion is now a Cauchy-Schwartz inequality.
Arguing as above, we find that
 \begin{eqnarray*}
     \displaystyle
\vertiii{A\bullet (M-\overline{M})\bullet B
}_{t,2n}^2
 &\leq& 
\int_0^t \EE\left[
\Vert A_s\Vert^{2n}\Vert B_s\Vert^{2n}~
\tr\left[
\partial_s\langle M-\overline{M}\,\vert\,\otimes\,\vert \,M-\overline{M}\,\rangle^{\,\sharp}_s
\right]^n
\right]^{1/n}~ds
    \end{eqnarray*}
This ends the proof of the theorem.
\cqfd

\section{Proof of the main results} \label{proof-section}

\subsection{Proof of theorem~\ref{theo-intro-phi-MM-Taylor}}\label{proof-theo-intro-phi-MM-Taylor}
We let $ \partial^{n}\phi_t(Q)$ be the collection of processes defined 
sequentially by the formulae
\begin{equation}\label{induction-partial-phi}
\begin{array}{l}
\displaystyle \partial^{n}\phi_t(Q)
\displaystyle=
\int_0^t
~\nabla\Lambda(\phi_s(Q))\cdot\partial^{n}\phi_s(Q)
~ds\\
\\
\displaystyle+\frac{1}{2}\sum_{1\leq k<n}\frac{n!}{k!(n-k)!}~\int_0^t
~\nabla^{2}\Lambda(\phi_s(Q))\cdot\left(\partial^{k}\phi_s(Q),\partial^{n-k}\phi_s(Q)\right)
~ds+n~\partial^{n-1}\MM_t(Q)
\end{array}\end{equation}
with  the collection of martingales
\begin{eqnarray*}
\displaystyle\partial^n\MM_t(Q)
&:=&\sum_{p+q=n}\sum_{(a,b)\in (\Ia_{p}\times\Ia_q)}~\frac{n!}{l(a)!l(b)!}
~\frac{1}{a!b!}~\MM^{(a,b)}_t(Q)
\\
\MM^{(a,b)}_t(Q)&:=&\int_0^t~\left(\left[\nabla^{l(a)}\varphi(\phi_s(Q))\cdot
\partial^{a}\phi_s(Q)\right]~d\Wa_s~
\left[\nabla^{l(b)}\Sigma_{\varphi}(\phi_s(Q))\cdot
\partial^{b}\phi_s(Q)\right]\right)_{ sym}
\end{eqnarray*}
Using (\ref{ref-nabla2-Lambda}) and Leibniz differential rule we find that
\begin{equation}\label{equivalent-form-derivatives}
\begin{array}{l}
\displaystyle \partial^{n}\phi_t(Q)
\displaystyle=\int_0^t\partial^n\Lambda(\phi_s(Q))~ds+n~\partial^{n-1}\MM_t(Q)
\end{array}
\end{equation}
Observe that
\begin{equation}\label{inductive-EE}
\begin{array}{l}
\displaystyle \partial^{n}\phi_t(Q)
=\displaystyle-\sum_{1\leq k<n}\frac{n!}{k!(n-k)!}~\int_0^t E_{s,t}(Q)
\left(\partial^{k}\phi_s(Q)\right)\left(\partial^{n-k}\phi_s(Q)\right)E_{s,t}(Q)^{\prime}
~ds\\
\displaystyle\hskip3cm
+n~\partial^{n-1}\MM_t(Q)+n\int_0^t~E_{s,t}(Q)~\left[\nabla\Lambda(\phi_s(Q))\cdot\partial^{n-1}\MM_s(Q)\right]~E_{s,t}(Q)^{\prime}~ds
\end{array}
\end{equation}
We check the estimates (\ref{estimate-partial-k})
 using an induction w.r.t. the parameter $n$.
 For $n=1$ we have
\begin{eqnarray*}
\displaystyle\partial\phi_t(Q)&=&\MM_t(Q)+\int_0^t~E_{s,t}(Q)~\left[\nabla\Lambda(\phi_s(Q))\cdot\MM_s(Q)\right]~E_{s,t}(Q)^{\prime}~ds
\end{eqnarray*}
 Using  (\ref{def-gamma}), (\ref{ref-M-bound-mn}) and (\ref{ref-M-bound-mn-2-integral}) we find that
 \begin{eqnarray*}
 \vertiii{ \partial\phi_t(Q)}_{2m}\vee  \left[
 t^{-1/2}\vertiii{ \partial\phi(Q)}_{t,2m}\right]&\leq& \vert Q\vert_+~G_{\nu}(\Vert Q\Vert,t)
 \end{eqnarray*}
We further assume $\partial^k\phi_t(Q)$ 
satisfy (\ref{estimate-partial-k})
for any integer $k\leq n$. In this situation, we have
$$
\begin{array}{l}
\displaystyle
\Vert\nabla^{l(a)}\varphi(\phi_s(Q))\cdot
\partial^{a}\phi_s(Q)\Vert
\displaystyle\leq  \vert Q\vert_-
\prod_{1\leq j\leq l(a)}\Vert
\partial^{a(j)}\phi_s(Q)\Vert
\end{array}
$$
Using corollary~\ref{cor-square-root-sigma-phi} we  check the estimate
$$
\Vert\nabla^{l(b)}\Sigma_{\varphi}(\phi_s(Q))\cdot
\partial^{b}\phi_s(Q)\Vert\leq \vert Q\vert_+ \prod_{1\leq j\leq l(b)}\Vert
\partial^{b(j)}\phi_s(Q)\Vert~
$$
On the other hand, combining (\ref{ref-M-bound-mn-2-integral}) with the induction hypothesis we have
$$
\begin{array}{l}
\vertiii{\partial^n \MM}_{t,2m}^2\leq\vert Q\vert~G_{\nu}(\Vert Q\Vert,t)~\left(1+t^{n+1}\right)
\Longrightarrow \vertiii{\partial^n \MM}_{t,2m}\leq (1+t)^{(n+1)/2}~\vert Q\vert~G_{\nu}(\Vert Q\Vert,t)
\end{array}
$$
Using Burkholder-Davis-Gundy inequality (\ref{lem-bdg-inequality}) and the induction hypothesis
we also check
$$
\vertiii{\int_0^t E_{s,t}(Q)d\left[\partial^{n-1}\MM_s(Q)\right]E_{s,t}(Q)^{\prime}}_{2m}\leq\vert Q\vert~G_{\nu}(\Vert Q\Vert,t)
$$
After some elementary manipulation we check that $\partial^{n+1}\phi(Q)$ 
satisfy (\ref{estimate-partial-k}) using the inductive formulae (\ref{induction-partial-phi}) and (\ref{inductive-EE}).
This ends the proof of (\ref{estimate-partial-k}).

Our next objective is to check that
the stochastic processes $\partial^n\phi_t(Q)$ and $\partial^n\MM_t(Q)$ defined in (\ref{induction-partial-phi}) coincide the $n$-th derivatives of the stochastic semigroup $\phi^{\epsilon}_t(Q)$ and the martingale $\MM_t^{\epsilon}(Q)$
at $\epsilon=0$. In addition, for $m,n\geq 1$ there exists some $\epsilon_{m,n}$ s.t. for  any $\epsilon\leq \epsilon_{m,n}$ we have
$$
\vertiii{\overline{\partial}^{\,n}\MM^{\epsilon}(Q)}_{t,m}\vee\vertiii{\overline{\partial}^{\,n}\phi^{\epsilon}(Q)}_{t,m}\leq \epsilon^n~\overline{e}(t)~\vert Q\vert
$$

We assume that the result has been proved up to rank $(n-1)$.
We have
$$
\begin{array}{l}
\displaystyle\phi^{\epsilon}_t(Q)
\displaystyle=Q+\int_0^t\Lambda\left(\phi_s(Q)+\overline{\partial}\phi^{\epsilon}_s(Q)\right)~ds\\
\displaystyle\hskip3cm+\epsilon~\int_0^t~\left(\varphi\left(\phi_s(Q)+\overline{\partial}\phi^{\epsilon}_s(Q)\right)~d\Wa_s~
\Sigma_{\varphi}(\phi_s(Q)+\overline{\partial}\phi^{\epsilon}_s(Q))\right)_{ sym}
\end{array}
$$
This implies that
$$
\begin{array}{l}
\displaystyle\overline{\partial}\phi^{\epsilon}_t(Q)
\displaystyle=\int_0^t\nabla\Lambda\left(\phi_s(Q)\right)\cdot\overline{\partial}\phi^{\epsilon}_s(Q\,ds+2^{-1}\int_0^t\nabla^2\Lambda\left(\phi_s(Q)\right)\cdot \overline{\partial}\phi^{\epsilon}_s(Q)\,ds\\
\qquad\qquad\displaystyle+\,\epsilon\int_0^t\,
\left(
\varphi\left(\phi_s(Q)\right)+
\overline{\nabla}\varphi\left[\phi_s(Q),\overline{\partial}\phi^{\epsilon}_s(Q)\right]
\right)\,d\Wa_s\,\left(
\Sigma_{\varphi}(\phi_s(Q))+\overline{\nabla}\Sigma_{\varphi}[\phi_s(Q),\overline{\partial}\phi^{\epsilon}_s(Q)]\right)_{ sym}
\end{array}
$$

On the other hand, we have
$$
\begin{array}{l}
\displaystyle\left(
\overline{\nabla}\varphi\left[\phi_s(Q),\overline{\partial}\phi^{\epsilon}_s(Q)\right]
~d\Wa_s~
\overline{\nabla}\Sigma_{\varphi}[\phi_s(Q),\overline{\partial}\phi^{\epsilon}_s(Q)]\right)_{ sym}\\
\\
\displaystyle=\sum_{1\leq k,l\leq n}\frac{1}{k!l!}
\displaystyle\left(\left[
\nabla^k\varphi(\phi_s(Q))\cdot\overline{\partial}\phi^{\epsilon}_s(Q)
\right]
~d\Wa_s~
\nabla^l\Sigma_{\varphi}(\phi_s(Q))\cdot\overline{\partial}\phi^{\epsilon}_s(Q)\right)_{ sym}+d\HH^{n,0}_s
\end{array}
$$
with the martingale
$$
\begin{array}{l}
\displaystyle d\HH^{n,0}_s
=\left(\left[
\overline{\nabla}\varphi[\phi_s(Q),\overline{\partial}\phi^{\epsilon}_s(Q)]
\right]
~d\Wa_s~
\overline{\nabla}^{n+1}\Sigma_{\varphi}[\phi_s(Q),\overline{\partial}\phi^{\epsilon}_s(Q)]\right)_{ sym}\\
\hskip4cm\displaystyle+
\displaystyle\left(\left[
\overline{\nabla}^{n+1}\varphi[\phi_s(Q),\overline{\partial}\phi^{\epsilon}_s(Q)]
\right]
~d\Wa_s~
\overline{\nabla}\Sigma_{\varphi}[\phi_s(Q),\overline{\partial}\phi^{\epsilon}_s(Q)]\right)_{ sym}
\end{array}
$$
Combining the first order estimate (\ref{Frob-Q-Q-estimates}) and the martingale estimates (\ref{ref-M-bound-mn})
with corollary~\ref{cor-square-root-sigma-phi} and corollary~\ref{cor-square-root-phi} we check that
$$
\vertiii{\HH^{n,0}}_{t,m}\leq \epsilon^n~\overline{e}(t)~\vert Q\vert
$$
We use the decomposition
$$
\begin{array}{l}
\displaystyle\sum_{1\leq k,l\leq n}\frac{1}{k!l!}\left(\left[
\nabla^k\varphi(\phi_s(Q))\cdot\partial\phi^{\epsilon}_s(Q)
\right]
~d\Wa_s~
\left[\nabla^l\Sigma_{\varphi}(\phi_s(Q))\cdot\partial\phi^{\epsilon}_s(Q)\right]\right)_{ sym}\\
\\
\displaystyle=\sum_{
\stackrel{1\leq p,q}{2\leq p+q<n}}~\frac{\epsilon^{p+q}}{(p+q)!}~
\sum_{(a,b)\in (\Ia_p\times\Ia_q)}
\frac{(p+q)!}{l(a)!l(b)!a!b!}~d\MM^{(a,b)}_s+d\HH^{n,1}_s\end{array}
$$
with the collection of martingales
$$
\begin{array}{l}
\displaystyle
d\HH^{n,1}_s\\
\\
\displaystyle:=\sum_{\stackrel{1\leq p,q}{n\leq  p+q\leq n^2}}
\sum_{(a,b)\in (\Ia_p\times\Ia_q)}
\left(\left[
\nabla^{l(a)}\varphi(\phi_s(Q))\cdot\partial^{(n,a)}\phi^{\epsilon}_s(Q)
\right]
~d\Wa_s~\left[
\nabla^{l(b)}\Sigma_{\varphi}(\phi_s(Q))\cdot\partial^{(n,b)}\phi^{\epsilon}_s(Q)\right]\right)_{ sym}
\end{array}
$$
defined in terms of functionals
\begin{eqnarray*}
\partial^{(n,a)}\phi^{\epsilon}_s(Q)&:=&\left(\partial^{(n,a(1))}\phi^{\epsilon}_s(Q),\ldots,\partial^{(n,a(p))}\phi^{\epsilon}_s(Q)\right)
\\
\mbox{\rm with}\quad\partial^{(n,a(i))}\phi^{\epsilon}_s(Q)&:=&1_{a(i)\not=n}~\frac{\epsilon^{a(i)}}{a(i)!}~\partial^{a(i)}\phi_s(Q)+
1_{a(i)=n}~\overline{\partial}^{n}\phi^{\epsilon}_s(Q)
\end{eqnarray*}
Observe that
$$
\begin{array}{l}
\displaystyle
\left(\overline{\nabla}[\phi_s(Q),\overline{\partial}\phi^{\epsilon}_s(Q)]~d\Wa_s~
\Sigma_{\varphi}(\phi_s(Q))\right)_{ sym}+\left(\varphi(\phi_s(Q))~d\Wa_s~
\overline{\nabla}\Sigma_{\varphi}[\phi_s(Q),\overline{\partial}\phi^{\epsilon}_s(Q)]\right)_{ sym}\\
\\
\displaystyle=\sum_{1\leq p< n}\frac{\epsilon^p}{p!}
\sum_{a\in \Ia_p}
\frac{p!}{l(a)!a!}~d\MM^{(a,0)}_s+\sum_{1\leq q< n}\frac{\epsilon^q}{q!}
\sum_{b\in \Ia_q}
\frac{q!}{l(b)!b!}~d\MM^{(0,b)}_s+d\HH^{n,2}_s
\end{array}
$$
with the collection of martingales
$$
\begin{array}{l}
\displaystyle
d\HH^{n,2}_s
\displaystyle:=\sum_{n\leq  p\leq n^2}
\sum_{a\in \Ia_p}
\left(\left[
\nabla^{l(a)}\varphi(\phi_s(Q))\cdot\partial^{(n,a)}\phi^{\epsilon}_s(Q)
\right]
~d\Wa_s~\Sigma_{\varphi}(\phi_s(Q))\right)_{ sym}\\
\\
\hskip3cm\displaystyle+\sum_{n\leq  q\leq n^2}
\sum_{b\in \Ia_q}
\left(\varphi(\phi_s(Q))
~d\Wa_s~~\left[
\nabla^{l(b)}\Sigma_{\varphi}(\phi_s(Q))\cdot\partial^{(n,b)}\phi^{\epsilon}_s(Q)\right]\right)_{ sym}
\end{array}
$$
Arguing as above and using the induction hypothesis we check that
$$
\forall i=1,2\qquad\vertiii{\HH^{n,i}}_{t,m}\leq \epsilon^n~\overline{e}(t)~\vert Q\vert$$
Finally observe that
$$
\begin{array}{l}
\displaystyle \nabla\Lambda\left(\phi_s(Q)\right)\cdot\overline{\partial}\phi^{\epsilon}_s(Q)+2^{-1}\nabla^2\Lambda\left(\phi_s(Q)\right)\cdot \overline{\partial}\phi^{\epsilon}_s(Q)\\
\\
\displaystyle =\sum_{1\leq k< n}\frac{\epsilon^{k}}{k!}~\nabla\Lambda\left(\phi_s(Q)\right)\cdot{\partial}^k\phi_s(Q)+\nabla\Lambda\left(\phi_s(Q)\right)\cdot\overline{\partial}^{n}\phi^{\epsilon}_s(Q)\\
\\
\displaystyle\hskip1cm+2^{-1}\sum_{2\leq k\leq n}\frac{\epsilon^{k}}{k!}~\sum_{1\leq p<k}~\frac{k!}{p!(k-p)!}~\nabla^2\Lambda\left(\phi_s(Q)\right)\cdot \left(\partial^{p}
\phi_s(Q),\partial^{k-p}
\phi_s(Q)\right)+H^{n}_s
\end{array}
$$
with the process
$$
\begin{array}{l}
\displaystyle
2H^{n}_s:=
\sum_{1\leq p,q,~n< p+q\leq 2n}\nabla^2\Lambda\left(\phi_s(Q)\right)\cdot \left(\partial^{(n,p)}
\phi^{\epsilon}_s(Q),\partial^{(n,q)}
\phi^{\epsilon}_s(Q)\right)
\Longrightarrow
\vertiii{H^n}_{t,m}\leq \epsilon^{n+1}~\overline{e}(t)~\vert Q\vert
\end{array}
$$
Using the above decompositions we check that
$$
\begin{array}{l}
\displaystyle\overline{\partial}\phi^{\epsilon}_t(Q)=\sum_{1\leq k< n}\frac{\epsilon^{k}}{k!}~
\int_0^t\nabla\Lambda\left(\phi_s(Q)\right)\cdot{\partial}^k\phi_s(Q)~ds+\int_0^t\nabla\Lambda\left(\phi_s(Q)\right)\cdot\overline{\partial}^{n}\phi^{\epsilon}_s(Q)~ds\\
\\
\displaystyle\hskip2cm+2^{-1}\sum_{2\leq k\leq n}\frac{\epsilon^{k}}{k!}~\sum_{1\leq p<k}~\frac{k!}{p!(k-p)!}~\int_0^t\nabla^2\Lambda\left(\phi_s(Q)\right)\cdot \left(\partial^{p}
\phi_s(Q),\partial^{k-p}
\phi_s(Q)\right)~ds
\\
\\
\displaystyle\hskip3cm+\sum_{2\leq k\leq n}~\frac{\epsilon^{k}}{k!}~k~\partial^{k-1}\MM_t(Q)+\int_0^t~H_s^n~ds+\epsilon~\left(\HH^{n,0}_t+\HH^{n,1}_t+\HH^{n,2}_t\right)
\end{array}
$$
Assuming the result is true at rank $(n-1)$ we check that
$$
\begin{array}{l}
\displaystyle\overline{\partial}^n\phi^{\epsilon}_t(Q)
\displaystyle=\int_0^t\nabla\Lambda\left(\phi_s(Q)\right)\cdot\overline{\partial}^{n}\phi^{\epsilon}_s(Q)~ds\\
\\
\hskip3cm\displaystyle+2^{-1}\frac{\epsilon^{n}}{n!}~\sum_{1\leq p<n}~\frac{n!}{p!(n-p)!}~\int_0^t\nabla^2\Lambda\left(\phi_s(Q)\right)\cdot \left(\partial^{p}
\phi_s(Q),\partial^{n-p}
\phi_s(Q)\right)~ds
\\
\\
\hskip4cm\displaystyle+\frac{\epsilon^{n}}{n!}~n~\partial^{n-1}\MM_t(Q)+\int_0^t~H_s^n~ds+\epsilon
~\left(\HH^{n,0}_t+\HH^{n,1}_t+\HH^{n,2}_t\right)
\end{array}
$$
This leads to the formula
$$
\begin{array}{l}
\displaystyle \frac{n!}{\epsilon^{n}}~\overline{\partial}^n\phi^{\epsilon}_t(Q)-\partial^n\phi_t(Q)\\
\\
\displaystyle=\int_0^t\nabla\Lambda\left(\phi_s(Q)\right)\cdot\left(
\frac{n!}{\epsilon^{n}}~\overline{\partial}^n\phi^{\epsilon}_t(Q)-\partial^n\phi_s(Q)\right)~ds+\frac{n!}{\epsilon^{n}}~\int_0^t~H_s^n~ds+\epsilon~\frac{n!}{\epsilon^{n}}~\left(\HH^{n,0}_t+\HH^{n,1}_t+\HH^{n,2}_t\right)
\end{array}
$$
The ends of the proof follows standard manipulations using Grownwall's lemma, thus it is skipped.
This ends the proof of the theorem.
\cqfd

\noindent Now we come to the proof of the bias estimate (\ref{intro-first-bias}).\label{intro-first-bias-proof}
For $n=2$ the stochastic process $\partial^{2}\phi_t(Q)$
 is defined by a diffusion equation (\ref{def-partial-phi}) that only depends on
the flows $(\phi_s(Q),\partial\phi_s(Q))$. 

Using (\ref{coupligng-eq-W}) we check that
$$
\partial_t\,
     \langle \MM(Q)\,\vert\,\otimes
\,\vert\, \partial\MM(Q)\rangle_t^{\,\sharp}=\phi_t(Q)\frownotimes\partial\Sigma(\phi_t(Q))+
\partial\phi_t(Q)\frownotimes\Sigma(\phi_t(Q))=\partial\left(\phi_t(Q)\frownotimes\Sigma(\phi_t(Q))\right)
$$
This yields the commutation formula
$$
   \langle \MM(Q)\,\vert\,\otimes
\,\vert\, \partial\MM(Q)\rangle=\partial \, \langle \MM(Q)\,\vert\,\otimes
\,\vert\, \MM(Q)\rangle\quad\mbox{\rm and}\quad \EE\left[   \langle \MM(Q)\,\vert\,\otimes
\,\vert\, \partial\MM(Q)\rangle\right]=0
$$
Using (\ref{inductive-EE}) and (\ref{prop-ra-lg-ll-gg-sharp}), this implies that
\begin{eqnarray*}
\displaystyle
\EE\left[\partial\phi_t(Q)\partial^2\phi_t(Q)\right]&=&0~=~
\EE\left[\partial\phi_t(Q)\right]
\\
\EE\left(\partial^{2}\phi_t(Q)\right)&=&-\int_0^t~E_{s,t}(Q)\EE\left[(\partial\phi_s(Q))^2\right]E_{s,t}(Q)^{\prime}~ds 
\Longrightarrow~~\EE(\partial^{3}\phi_t(Q))=0
\end{eqnarray*}
This ends the proof of the bias estimate (\ref{intro-first-bias}).\cqfd

\subsection{Proof of theorem~\ref{th1-estimates-Ln}}\label{th1-estimates-Ln-proof}

\begin{lem}[Gronwall-Riccati lemma]\label{gronwall-riccati}
Let $f_t$ be some function satisfying differential 
inequality
\begin{equation}\label{f13}
\partial_tf_t\leq  a~f_t-  b~f_t^2+ c
\end{equation}
for some $ a\in \RR$, and $ b,c>0$. We set 
$
h:=(2 b)^{-1}( a+\sqrt{ a^2+4 b c})
$. 
 In this situation, we have the following estimates:
 If $f_0=h$ then $f_t\leq h$. When 
 $f_0\not=h$ we have
$
0\leq f_t\leq h+\left(
f_0- h\right)_+~e^{-t\sqrt{ a^2+4 b c}}$.
In addition, when $f_0>h$ for any $t>0 $ we have the uniform estimate
\begin{eqnarray*}
f_t< h+\left( b^{-1}{\sqrt{ a^2+4 b c}}~\left[\exp{\left(t\sqrt{ a^2+4 b c}\right)}-1\right]^{-1}\wedge \vert f_0-h\vert \right)
\end{eqnarray*}

\end{lem}
\proof
We let $g_t$ be the non negative solution of the Riccati equation
$$
\partial_tg_t=  a~g_t-  b~g_t^2+ c
$$
When $ b=0$, by Gronwall's lemma we have $f_t\leq g_t$.
We further assume that $ b\not=0$.
Notice that
\begin{eqnarray*}
\partial_t(f_t-g_t)&\leq& (f_t-g_t)\left[\left( a-2 b g_t\right)- b (f_t-g_t)\right]\leq  \left( a-2 b g_t\right)~(f_t-g_t)
\end{eqnarray*}
By Gronwall's lemma this implies that $f_t\leq g_t$ for any $t\geq 0$, as soon as $f_0=g_0$.
The estimates are now a direct consequence of the properties of the one dimensional Riccati equation.
This ends the proof of the lemma.
\cqfd

Now we come to the proof of theorem~\ref{th1-estimates-Ln}.

{\bf Proof of theorem~\ref{th1-estimates-Ln}:} Taking the trace and using (\ref{trQ2}), for any $Q\in \Sa_r^0$ we find that
\begin{eqnarray*}
\tr\left(\Lambda(Q)\right)
&\leq &\tr(R)+2\rho(A)~\tr\left(Q\right)-r^{-1}\tr\left(Q\right)^2
\end{eqnarray*}
This shows that
$$
(p_t,q_t):=(\tr(P_t),\tr(Q_t))\quad\mbox{\rm with}\quad
(P_t,Q_t):=\left(\phi_t(P),\phi^{\epsilon}_t(P)\right)
$$
satisfy the differential inequalities
$$
\partial_tp_t\leq \gamma(p_t):=a~p_t-b~p_t^2+c
\quad\mbox{\rm and}\quad
dq_t\leq \gamma(q_t)~dt+\epsilon~dm^{\epsilon}_t
$$
with $m^{\epsilon}_t:=\tr(\MM^{\epsilon}_t)$ and  $(a,b,c)=(2\rho(A),1/r,\tr(R))$.
Also notice that
$$
\partial_t\langle m^{\epsilon}\,\vert\, m^{\epsilon}\rangle_t=\tr\left(Q_t\Sigma(Q_t)\right)\leq 
\tau(q_t):=q_t(\lambda_{max}(R)+q_t^2)
$$
In this notation, for any $n\geq 1$ we have
\begin{eqnarray*}
dq_t^{\,n}
&\leq & n~\left[a~q_t^{\,n}-b_n q_t^{\,n+1}+c_n~q_t^{\,n-1}\right]~dt+n~q_t^{\,n-1}~\epsilon~dm^{\epsilon}_t
\\
\mbox{\rm with}\quad
b_n&:=&b-\frac{n-1}{2}~\epsilon^2\quad\mbox{\rm and}\quad c_n:=c+\frac{n-1}{2}~\epsilon^2~\lambda_{max}(R)
\end{eqnarray*}
Taking expectation we find  the Riccati differential inequality
$$
\displaystyle f_t:=
\EE\left[q_t^n\right]^{\frac{1}{n}}\Longrightarrow\partial_tf_t\leq 
 a~f_t- b_n~f_t^2+ c_n
$$
We further assume that $\epsilon$ is chosen so that
$
r~(n-1)~\epsilon^2<2
$. Using lemma~\ref{gronwall-riccati} we have
$$
f_t~\leq~ b_n^{-1}( a+\sqrt{ a^2+4 b_n c_n})+f_0 ~\leq~ 2\left(\frac{a}{b_n}+\sqrt{\frac{c_n}{b_n}}\right)+f_0
$$
This ends the proof of  the r.h.s. estimate stated in (\ref{trace-estimates-unif}).
We also have the estimates
$$
\begin{array}{l}
\displaystyle \overline{q}_t:=\sup_{s\in [0,t]}{q_s}\leq  ct+q_0+\int_0^t~a_+~\overline{q}_s~~ds+\epsilon~\overline{m}_t\quad\mbox{\rm with}\quad \overline{m}_t:=\sup_{s\in [0,t]}{m^{\epsilon}_s}
\\
\\
\displaystyle\Longrightarrow\EE\left(\overline{q}_t^n\right)^{1/n}\leq  
ct+\EE\left({q_0^{n}}\right)^{1/n}+\int_0^t~a_+~
\EE\left(\overline{q}_s^n\right)^{1/n}~ds+\epsilon~\EE\left(\overline{m}_t^n\right)^{1/n}
\end{array}
$$
Combing Burkholder-Davis-Gundy inequality with  the r.h.s. estimate stated in (\ref{trace-estimates-unif}) we find that
\begin{eqnarray*}
r~(3n-1)~\epsilon^2<2&\Longrightarrow&
\EE\left(\overline{m}_t^{2n}\right)^{1/n}\leq c_1
\int_0^t \EE(q^n_s)^{1/n}ds+c_2
\int_0^t \EE(q^{3n}_s)^{1/n}ds\leq ~t~\vert P\vert_+\\
&\Longrightarrow&\displaystyle\EE\left(\overline{q}_t^n\right)^{1/n}\leq 
\EE\left({q_0^{n}}\right)^{1/n}+\int_0^t~a_+~
\EE\left(\overline{q}_s^n\right)^{1/n}~ds+~t~\vert P\vert_+
\end{eqnarray*}
The l.h.s. estimate stated in (\ref{trace-estimates-unif}) is now a direct consequence of Gronwall's.
Now we come to the proof of (\ref{Frob-Q-Q-estimates}). We set 
$$
\widehat{\partial}\,\phi^{\epsilon}_t(P)=\epsilon^{-1}~\overline{\partial}\phi^{\epsilon}_t(P)=\epsilon^{-1}\left[\phi^{\epsilon}_t(P)-\phi_t(P)\right]\Longrightarrow
\phi_s(P)+\frac{\epsilon}{2}~\widehat{\partial}\,\phi^{\epsilon}_s(P)=\frac{\phi_u(P)+\phi^{\epsilon}_u(P)}{2}.
$$
The mean value formula (\ref{def-Lambda}) also yields
\begin{eqnarray*}
\widehat{\partial}\,\phi^{\epsilon}_t(P)
&=&\int_0^t E^{\epsilon}_{s,t}(P)~d\MM^{\epsilon}_s(P)~E^{\epsilon}_{s,t}(P)^{\prime}
\end{eqnarray*}
Using the integration by part we find that
$$
\widehat{\partial}\,\phi^{\epsilon}_t(P)
\displaystyle=\MM^{\epsilon}_t(P)+
\int_0^t~E^{\epsilon}_{s,t}(P)\left[\nabla\Lambda\left(\phi_s(P)+\frac{\epsilon}{2}~\widehat{\partial}\,\phi^{\epsilon}_s(P)\right)\cdot\MM^{\epsilon}_s(P)\right]E^{\epsilon}_{s,t}(P)^{\prime}~ds
$$
Using the estimates  (\ref{trace-estimates-unif}) we  check  the almost sure estimate
$$
\displaystyle
\Vert \widehat{\partial}\,\phi^{\epsilon}(P)\Vert_t\leq   \overline{e}(t)~\Vert\MM^{\epsilon}\Vert_t~\vert P\vert_+
$$
The end of the proof of the l.h.s. estimate in (\ref{Frob-Q-Q-estimates}) is now a consequence of the estimates
(\ref{trace-estimates-unif}) and (\ref{ref-M-bound-mn}). Using (\ref{trace-estimates-unif}) and (\ref{ref-M-bound-mn-2-integral}) we also have
$$
\Vert \widehat{\partial}\,\phi^{\epsilon}_t(P)\Vert^2_{2n}\leq c_1~
\int_0^t~e^{4\rho(A)(t-s)}~
\left(1+\EE\left(\Vert \phi_s^{\epsilon}(P)\Vert^{6n}\right)^{1/n}\right)~ds\leq c_1~\frac{1-e^{4\rho(A)t}}{4\rho(A)}
~\vert P\vert_+
$$
as soon as  $(18n-1)~\epsilon^2<2/r$.
This ends the proof of the theorem.\cqfd

\subsection{Proof of theorem~\ref{th2}}\label{th2-proof}

Using (\ref{monotone-prop}) and  (\ref{ref-M-bound-mn})  we have
$
\vertiii{\MM(P)}_{t,n}\leq c~\sqrt{t}~\vert P\vert_+
$.
On the other hand, we have
$$
\begin{array}{l}
\displaystyle\Lambda(\phi^{\epsilon}_t(P))-\Lambda(\phi_t(P))-\epsilon~\nabla\Lambda(\phi_t(P))\cdot \partial\phi_t(P)
=\nabla\Lambda(\phi_t(P))\cdot 
 \overline{\partial}^2\phi^{\epsilon}_t(P)-\left(\overline{\partial}\phi^{\epsilon}_t(P)\right)^2
\end{array}
$$
This yields the decomposition
$$
\begin{array}{l}
\displaystyle \overline{\partial}^2\phi_t^{\epsilon}(P)
\displaystyle=\int_0^t~\nabla\Lambda(\phi_s(P))\cdot 
 \overline{\partial}^2\phi^{\epsilon}_s(P)~ds-\int_0^t~
\left( \overline{\partial}\phi_s^{\epsilon}(P)\right)^2~ds
\displaystyle+\epsilon~\overline{\partial}\,\MM_t^{\epsilon}(P)\end{array}
$$
from which we check that
$$
\begin{array}{l}
\displaystyle \overline{\partial}^2\phi_t^{\epsilon}(P)
\displaystyle=
-\int_0^t~E_{s,t}(P)~\left(\overline{\partial}\phi_s^{\epsilon}(P)\right)^2~E_{s,t}(P)^{\prime}~ds
\displaystyle+\epsilon~\int_0^tE_{s,t}(P)~d\left[\overline{\partial}\,\MM_s^{\epsilon}(P)\right]~E_{s,t}(P)^{\prime}
\end{array}
$$
This implies that\label{proof-order-1-sobolev}
$$
\begin{array}{l}
\displaystyle \EE\left(\overline{\partial}^2\phi_t^{\epsilon}(P)\right)=\EE\left(\phi^{\epsilon}_t(P)-\phi_t(P)\right)
= -\int_0^t~E_{s,t}(P)~\EE\left(\left[\overline{\partial}\phi_s^{\epsilon}(P)\right]^2\right)~E_{s,t}(P)^{\prime}~ds\leq 0\\
\\
\displaystyle\Longrightarrow \EE\left(\phi^{\epsilon}_t(P)\right)=\phi_t(P)-\int_0^t~E_{s,t}(P)~\EE\left(\left[\overline{\partial}\phi_s^{\epsilon}(P)\right]^2\right)~E_{s,t}(P)^{\prime}~ds\leq \phi_t(P)
\end{array}
$$

This ends the proof of the l.h.s. estimate in (\ref{order-1-sobolev}). 
We check by an integration by part that 
$$
\begin{array}{l}
\displaystyle \overline{\partial}^2\phi_t^{\epsilon}(P)=
\epsilon~\overline{\partial}\,\MM_t^{\epsilon}(P)\\
\\
\displaystyle+\epsilon~\int_0^tE_{s,t}(P)\left(\nabla\Lambda(\phi_s(P))\cdot \overline{\partial}\,\MM_s^{\epsilon}(P)\right)E_{s,t}(P)^{\prime}~ds-\int_0^t~E_{s,t}(P)~\left(\overline{\partial}\phi_s^{\epsilon}(P)\right)^2~E_{s,t}(P)^{\prime}~ds
\end{array}
$$
Using (\ref{monotone-prop}) and the estimates (\ref{def-gamma}) we check the almost sure estimate
$$
\Vert\overline{\partial}^2\phi^{\epsilon}(P)\Vert_t
\leq \epsilon~G_{\nu}(\Vert P\Vert,t)~
\Vert \overline{\partial}\MM^{\epsilon}(P) \Vert_t~\Vert P\Vert+G_{\nu}(\Vert P\Vert,t)~
\Vert\overline{\partial}\phi^{\epsilon}(P)\Vert_t^2
$$
On the other hand, using (\ref{pr-fcv}) we have
$$
\begin{array}{l}
\displaystyle\vertiii{\overline{\partial}\,\MM^{\epsilon}(P)}_{t,n}^2
\displaystyle\leq 
\int_0^t~\chi^{\epsilon}_n(s)~\left[\vertiii{
\phi^{\epsilon}_s(P)-\phi_s(P)}_n+\vertiii{
\phi^{\epsilon}_s(P)^2-\phi_s(P)^2}_n\right]~ds
\end{array}
$$
with the function
$$
\begin{array}{l}
\displaystyle
\chi^{\epsilon}_n(s)
:=1+\Vert P\Vert^2+\vertiii{ \phi^{\epsilon}_s(P)}_{2n}^2
\end{array}
$$
Using (\ref{trace-estimates-unif}) and (\ref{Frob-Q-Q-estimates}) we have
$$
\begin{array}{l}
\displaystyle\vertiii{\overline{\partial}\,\MM^{\epsilon}(P)}_{t,n}^2
\displaystyle\leq ~t~
\left[\vert P\vert_++\vertiii{ \phi^{\epsilon}(P)}_{t,2n}^3\right]
\vertiii{\overline{\partial}
\phi^{\epsilon}(P)}_{t,2n}\leq \epsilon~\overline{e}(t)~\vert P\vert_+
\end{array}
$$
as soon as $\epsilon\leq \epsilon_n$, for some $\epsilon_n$.
The end of the proof of the r.h.s. estimate in (\ref{order-1-sobolev}) is now a consequence of (\ref{Frob-Q-Q-estimates}).
By proposition 2.2 in~\cite{aps-2016}, for any $P\in\Sa_{r}^0$ and any $0\leq s\leq t$ we have
\begin{equation}\label{la-ref}
\partial_{s}\phi_{s,t}(P)=-\Lambda(\phi_{s,t}(P))=-\nabla\phi_{s,t}(P)\cdot \Lambda(P)
\end{equation}
Using theorem~\ref{lem-frechet-derivatives}  and applying the Ito formula (\ref{ito-formula})  to the function
$\Upsilon(P)=\phi_{s,t}(P)$ we also have
$$
\begin{array}{rcl}
\displaystyle\nabla\Upsilon(P)\cdot H&=&E_{t-s}(P)HE_{t-s}(P)^{\prime}\\
\displaystyle\frac{1}{2}\nabla^2\Upsilon(P)\cdot H&=&
-\left[E_{t-s}(P)H \Gamma_{t-s}(P)^{1/2}\right]\left[E_{t-s}(P)H \Gamma_{t-s}(P)^{1/2}\right]^{\prime}\\
\displaystyle-~(P\frownotimes \Sigma(P))\nabla^2\Upsilon(P)&=&2\left[(E_{t-s}(P)\otimes \Gamma_{t-s}(P)^{1/2})(P\frownotimes \Sigma(P))(\Gamma_{t-s}(P)^{1/2}\otimes E_{t-s}(P) )^{\prime}\right]^{\flat}\\
&=&\displaystyle
\Omega_{t-s}(P)\end{array}$$
This implies that
\begin{equation}\label{ref-tcl}
\begin{array}{l}
\displaystyle
\phi^{\epsilon}_t(P)=\phi_{t}\left(P\right)
+\epsilon~\int_0^t~E_{t-s}(\phi^{\epsilon}_s(P))~ d\MM^{\epsilon}_{s}(P)~E_{t-s}(\phi^{\epsilon}_s(P))^{\prime}
-\frac{\epsilon^2}{2}~
\int_0^t\Omega_{t-s}\left(\phi^{\epsilon}_s(P)\right)~ds
\end{array}
\end{equation}
On the other hand  we have
\begin{eqnarray}\label{DL-2-E}
E_t(\phi^{\epsilon}_s(P))
&=&E_t(\phi_s(P))+\epsilon~\partial (E_t\circ\phi_s)(P)+\overline{\partial}^2(E_t\circ\phi_s)(P)
\end{eqnarray}
\begin{eqnarray*}
\mbox{\rm with the matrices}\quad\partial (E_t\circ\phi_s)(P)&:=&\nabla E_t(\phi_s(P))\cdot \partial\phi_s(P)\Rightarrow\EE\left(\partial (E_t\circ\phi_s)(P)\right)=0
\\
\overline{\partial}^2(E_t\circ\phi_s)(P)&:=&\nabla E_t(\phi_s(P))\cdot \overline{\partial}^2\phi^{\epsilon}_s(P)
+\overline{\nabla}^2E_t\left[\phi^{\epsilon}_s(P),\overline{\partial}\phi^{\epsilon}_s(P)\right]
\end{eqnarray*}
Using (\ref{nabla-En})
we check that
$$
\vertiii{ \overline{\nabla}E_t\left[\phi^{\epsilon}_s(P),\overline{\partial}\phi^{\epsilon}_s(P)\right]}\leq 
\overline{g}(t)~\quad\mbox{\rm and}\quad
\Vert \overline{\partial}^2E_t(\phi^{\epsilon}_s(P))\Vert\leq 
\overline{g}(t)~\left[\Vert \overline{\partial}^2\phi^{\epsilon}_s(P)
\Vert+~\Vert \overline{\partial}\phi^{\epsilon}_s(P)\Vert^2\right]
$$
In the same vein, using (\ref{estimate-nabla-Pi}) we also have
\begin{eqnarray*}
\Pi_t(\phi^{\epsilon}_s(P))
&=&\Pi_t(\phi_s(P))+\epsilon~\partial (\Pi_t\circ\phi_s)(P)+\overline{\partial}^2\Pi_t(\phi^{\epsilon}_s(P))
\end{eqnarray*}
with $\EE\left(\partial (\Pi_t\circ\phi_s)(P)\right)=0$ and the estimates
$$
\begin{array}{l}
\displaystyle\Vert \overline{\nabla}\Pi_t\left[\phi^{\epsilon}_s(P),\overline{\partial}\phi^{\epsilon}_s(P)\right]\Vert\leq G(t)~
\vert \phi^{\epsilon}_s(P)\vert_+~\vert \overline{\partial}\phi^{\epsilon}_s(P)\vert_+^2 \\
\\
\displaystyle\Vert \overline{\partial}^2\Pi_t(\phi^{\epsilon}_s(P))\Vert\leq 
G(t)~\left[\vert P\vert_++\vert \overline{\partial}\phi^{\epsilon}_s(P)\vert_++
\vert \phi^{\epsilon}_s(P)\vert_+\right]~\left[~\Vert \overline{\partial}^2\phi^{\epsilon}_s(P)\Vert+ \Vert \overline{\partial}\phi^{\epsilon}_s(P)\Vert^2\right]
\end{array}
$$
Combining (\ref{estimate-partial-k}) and (\ref{monotone-prop})
with (\ref{def-gamma}) and (\ref{estimate-nabla-Pi}) we check that
$$
\begin{array}{l}
\displaystyle
\partial (\Omega_{t}\circ\phi_s)(P)\\
\\\displaystyle
=E_t(\phi_s(P))~\partial\,(\Pi_t\circ\phi_s)(P)~E_t(\phi_s(P))^{\prime}-2
~\left(E_t(\phi_s(P))~\Pi_t(\phi_s(P))~\partial (E_t\circ\phi_s)(P)^{\prime}\right)_{ sym}\\
\\
\Longrightarrow\Vert \partial \Omega_{t}\left(\phi_s(P)\right)\Vert\leq 
\overline{g}_{\nu}(\Vert P\Vert,t)~\Vert P\Vert_{+}
\end{array}$$
We also have
$$
\begin{array}{l}
\displaystyle\Vert\overline{\partial}^2 \Omega_{t}\left(\phi^{\epsilon}_s(P)\right)\Vert
\displaystyle\leq \overline{g}(t)~\left[\vert P\vert_++\vert \overline{\partial}\phi^{\epsilon}_s(P)\vert_++
\Vert \phi^{\epsilon}_s(P)\Vert _+\right]~\left[~\Vert \overline{\partial}^2\phi^{\epsilon}_s(P)\Vert+ \Vert \overline{\partial}\phi^{\epsilon}_s(P)\Vert^2\right]
\end{array}
$$
Also notice that
\begin{eqnarray*}
\displaystyle\partial (\Omega_{t}\circ\phi_s)(P)
&=&\nabla\Omega_{t}(\phi_s(P))\cdot \partial\phi_s(P)\Longrightarrow\EE\left(\partial (\Omega_{t}\circ\phi_s)(P)\right)=0\\
\displaystyle\nabla\Omega_{t}(\phi_s(P))\cdot H\displaystyle
&:=&E_t(\phi_s(P))~\left[\nabla\Pi_t(\phi_s(P))\cdot H\right]~E_t(\phi_s(P))^{\prime}\\
&&\displaystyle\hskip3cm-2
~\left(E_t(\phi_s(P))~\Pi_t(\phi_s(P))~\left[\nabla E_t(\phi_s)(P))\cdot H\right]^{\prime}\right)_{ sym}
\end{eqnarray*}
We conclude that
$$
\EE\left[\phi^{\epsilon}_t(P)\right]-\phi_{t}\left(P\right)+\frac{\epsilon^2}{2}~
\int_0^t \Omega_{t-s}\left(\phi_s(P)\right)
~ds
=
-\frac{\epsilon^2}{2}~
\int_0^t\EE\left[\overline{\partial}^2 \Omega_{t-s}\left(\phi^{\epsilon}_s(P)\right)\right]~ds
$$
with
 the second order remainder
\begin{eqnarray*}
\displaystyle
\Vert\EE\left[\overline{\partial}^2 \Omega_{t-s}\left(\phi^{\epsilon}_s(P)\right)\right]\Vert
&\leq& \overline{g}(t-s)~\overline{e}(s)~\epsilon^{3/2}~\vert P\vert
\end{eqnarray*}
as soon as $\epsilon\leq \epsilon_0$, for some $\epsilon_0$. This yields the estimate
\begin{equation}\label{th2-estimate-tcl-bias}
\Vert\EE\left[\phi^{\epsilon}_t(P)\right]-\phi_{t}\left(P\right)+\frac{\epsilon^2}{2}~
\int_0^t \Omega_{t-s}\left(\phi_s(P)\right)\Vert
\displaystyle\leq \epsilon^{7/2}~\overline{e}(t)~~\vert P\vert
\end{equation}
We check (\ref{equivalent-bias-terms})
combining  (\ref{intro-first-bias}) with  (\ref{def-partial-phi}) and (\ref{th2-estimate-tcl-bias}). The proof of the theorem is complete.\cqfd

\subsection{Proof of theorem~\ref{tcl+initial-condition}}\label{tcl+initial-condition-proof}

The proof of the theorem is based on the following technical lemma.
\begin{lem}
The mapping $Q\mapsto \MM_t(Q)$ is smooth with first and second order derivatives
given by the martingales
$$
\begin{array}{rcl}
\displaystyle
\nabla\MM_t(Q)\cdot H&=&\displaystyle\int_0^t\left(
(\varphi\circ\phi_s)(Q)~d\Wa_s~\left[\nabla(\Sigma_{\varphi}\circ\phi_s)(Q)\cdot H\right]\right)_{ sym}\\
&&\hskip2cm \displaystyle+ \int_0^t\left(
\left[\nabla (\varphi\circ\phi_s)(Q)\cdot H\right]~d\Wa_s~(\Sigma_{\varphi}\circ\phi_s)(Q)\right)_{ sym}\\
&&\\
\displaystyle \nabla^2\MM_t(Q)\cdot H
&=&\displaystyle\int_0^t\left(
(\varphi\circ\phi_s)(Q)~d\Wa_s~\left[
\nabla^2(\Sigma_{\varphi}\circ\phi_s)(Q)\cdot H\right]\right)_{ sym}\\
\displaystyle&&\hskip1cm+\displaystyle \int_0^t~\left(
\left[\nabla^2 (\varphi\circ\phi_s)(Q)\cdot H\right]~d\Wa_s~
(\Sigma_{\varphi}\circ\phi_s)(Q)\right)_{ sym}\\
\displaystyle&&\hskip2cm\displaystyle +2~\int_0^t~\left(\left(
\nabla (\varphi\circ\phi_s)(Q)\cdot H\right)~d\Wa_s~\left(\nabla(\Sigma_{\varphi}\circ\phi_s)(Q)\cdot H
 \right)\right)_{ sym}
\end{array}
$$
In addition, for any $i=1,2$ we have the estimates 
\begin{equation}\label{estimates-M-nabla-M}
\vertiii{\MM(Q)}_{t,n}\leq \sqrt{t}~\vert Q\vert_+\qquad
\vertiii{\nabla^i\MM(Q)}_{t,n}\vee 
\left[\vert H\vert_+^{-1}\vertiii{\overline{\nabla}^3\MM\{Q,H\}}_{t,n}\right]
\leq~\overline{e}(t)~\vert Q\vert
\end{equation}

\end{lem}

\proof
By theorem~\ref{lem-square-root-taylor} and theorem~\ref{lem-frechet-derivatives}, the mapping $\varphi\circ\phi_t$ is smooth and we have
\begin{eqnarray*}
\nabla (\varphi\circ\phi_t)(Q)\cdot H&=&\nabla\varphi\left(\phi_t(Q)\right)\cdot(
\nabla\phi_t(Q)\cdot H)
\Longrightarrow\vertiii{\nabla (\varphi\circ\phi_t)(Q)}\leq e(t)~\vert Q\vert_-\\
\nabla^2 (\varphi\circ\phi_t)(Q)\cdot (H,H)
&=&\nabla^2\varphi\left(\phi_t(Q)\right)\cdot\left(
\nabla\phi_t(Q)\cdot H\right)+\nabla\varphi\left(\phi_t(Q)\right)\cdot\left(\nabla^2\phi_t(Q)\cdot H\right)\\
&&\Longrightarrow\vertiii{\nabla^2 (\varphi\circ\phi_t)(Q)}\leq  e(t)~\vert Q\vert_-
\end{eqnarray*}
The first and the second order remainder terms are given by
$$
\begin{array}{l}
\displaystyle\overline{\nabla}^2
 (\varphi\circ\phi_t)[Q,H]
\displaystyle =\nabla\varphi\left(\phi_t(Q)\right)\cdot
\overline{\nabla}^2\phi_t\left[Q,H\right]\\
\\
\displaystyle\hskip2cm+2^{-1}
\nabla^2\varphi\left(\phi_t(Q)\right)~\cdot\left(\overline{\nabla}\phi_t\left[Q,H\right],\overline{\nabla}\phi_t\left[Q,H\right]\right)+\overline{\nabla}^3\varphi\left[\phi_t(Q),\overline{\nabla}\phi_t\left[Q,H\right]\right]\\
\\
\Longrightarrow\vertiii{ \overline{\nabla}^2
 (\varphi\circ\phi_t)[Q,H]}
  \leq e(t)~\vert Q\vert_-~(1+\Vert H\Vert)
 \end{array}
$$
We also have
$$
\begin{array}{l}
\displaystyle\overline{\nabla}^3
 (\varphi\circ\phi_t)[Q,H]
=
\displaystyle 
\nabla^2\varphi\left(\phi_t(Q)\right)~\cdot\left(\overline{\nabla}^{2}\phi_t[Q, H],~
\left[\nabla\phi_t(Q)\cdot
H\right]\right)\\
\\
\hskip4cm\displaystyle+2^{-1}
\nabla^2\varphi\left(\phi_t(Q)\right)~\cdot\left(
\overline{\nabla}^{2}\phi_t[Q, H],~
\overline{\nabla}^{2}\phi_t[Q, H]\right)\\
\\
\hskip5cm\displaystyle+\overline{\nabla}^3\varphi\left[\phi_t(Q),\overline{\nabla}\phi_t\left[Q,H\right]\right]
+\nabla\varphi\left(\phi_t(Q)\right)\cdot
\overline{\nabla}^{3}\phi_t[Q, H]\\
\\
\Longrightarrow\vertiii{ \overline{\nabla}^3
 (\varphi\circ\phi_t)[Q,H]}\leq e(t)~\vert Q\vert_-~(1+\Vert H\Vert)~
 \end{array}
$$
In the same vein we prove that $\Sigma_{\varphi}\circ\phi_t$ is smooth 
with the first and second derivatives and the first and second order remainder terms are such that
$$
\forall k=1,2\qquad\vertiii{ \overline{\nabla}^k\left(\Sigma_{\varphi}\circ \phi_t\right)[Q,H]}\leq e(t)~\vert H\vert_+
\vert Q\vert_+
$$
Using (\ref{ref-M-bound-mn}) we prove the estimates (\ref{estimates-M-nabla-M}). This ends the proof of the lemma. \cqfd

We are now in position to prove theorem~\ref{tcl+initial-condition}.

{\bf Proof of theorem~\ref{tcl+initial-condition}:} The mapping $Q\mapsto \partial\phi_t(Q)$ is smooth with 
$$
\begin{array}{l}
\nabla \partial\phi_t(Q)\cdot H
=\displaystyle
\int_0^t~E_{s,t}(Q)~d\,\left(\nabla\MM_s(Q)\cdot H\right)~E_{s,t}(Q)^{\prime}\\
\displaystyle\hskip4cm+
2~\int_0^t~\left[E_{s,t}(Q)~d\MM_s(Q)~
\left(\nabla E_{s,t}(Q)\cdot H\right)^{\prime}\right]_{ sym}
\end{array}
$$
$$
\begin{array}{l}
\displaystyle\nabla^2 \partial\phi_t(Q)\cdot H
\displaystyle=\int_0^t~E_{s,t}(Q)~d\,\left(\nabla^2\MM_s(Q)\cdot H\right)~E_{s,t}(Q)^{\prime}\\
\displaystyle\hskip4cm+4~\int_0^t~\left(E_{s,t}(Q)~d\,\left(\nabla\MM_s(Q)\cdot H\right)
\left(\nabla E_{s,t}(Q)\cdot H\right)^{\prime}\right)_{ sym}\\
\displaystyle\hskip4.5cm+\displaystyle 2\int_0^t~\left(\nabla E_{s,t}(Q)\cdot H\right)~d\MM_s(Q)~\left(
\nabla E_{s,t}(Q)\cdot H\right)^{\prime} \\
\displaystyle\hskip5cm+2\int_0^t~\left[\left(\nabla^2 E_{s,t}(Q)\cdot H
\right)~d\MM_s(Q)~E_{s,t}(Q)^{\prime}\right]_{ sym}
\end{array}
$$
After some manipulations we also check that
$$
\forall k=1,2\qquad
\vertiii{\overline{\nabla}^k \partial\phi\{Q,H\}}_{t,m}
\displaystyle\leq~\overline{e}(t)~\vert Q\vert~
\vert H\vert_+
$$
We fix some matrix $H$ and we set
$
\psi^{\epsilon}(Q):=Q+\epsilon~H
$.
Using (\ref{order-1-sobolev})~and~(\ref{estimate-nabla-k-phi}) we have
\begin{eqnarray*}
\phi^{\epsilon}_t\left(\psi^{\epsilon}(Q)\right)&:=&\phi^{\epsilon}_t\left(Q+\epsilon~H\right)=\phi_t\left(Q\right)+\epsilon ~{\partial}(\phi_t\circ\psi)(Q)+\overline{\partial}^{2}(\phi^{\epsilon}_t\circ\psi^{\epsilon})(Q)
\end{eqnarray*}
with
\begin{eqnarray*}
{\partial}(\phi_t\circ\psi)(Q)&:=&\left[\nabla\phi_t(Q)\cdot H+\partial\phi_t\left(Q\right)\right]{\Longrightarrow}~\Vert {\partial}(\phi_t\circ\psi)(Q)\Vert\leq \overline{e}(t)~\vert Q\vert_+\\
\overline{\partial}^{2}(\phi^{\epsilon}_t\circ\psi^{\epsilon})(Q)&:=&\overline{\nabla}^2\phi_t[Q,\epsilon H]
+
\epsilon~\overline{\nabla}\partial\phi_t\left[Q,\epsilon~H\right]+\overline{\partial}^{2}\phi_t^{\epsilon}\left(Q+\epsilon~H\right)
\end{eqnarray*}
Using the r.h.s. estimate in (\ref{order-1-sobolev}) and (\ref{estimate-nabla-k-phi}) we also find that
$$
\begin{array}{l}
\vertiii{\overline{\partial}^{2}(\phi^{\epsilon}\circ\psi^{\epsilon})(Q)}_{t,n}
\leq\epsilon^{3/2}~\overline{e}(t)~\vert Q\vert~\vert H\vert_+\end{array}$$
This ends the proof of the estimate (\ref{estimate-tcl}) for $k=1$.
We also have
$$
\begin{array}{l}
\left[(\phi^{\epsilon}_t\circ\psi^{\epsilon})(Q)-\phi_t\left(Q\right)\right]^2
\displaystyle-\epsilon^2~\left[
 {\partial}(\phi_t\circ\psi)(Q)\right]^2\\
\\
\displaystyle=\epsilon ~\left[
{\partial}(\phi_t\circ\psi)(Q)~\overline{\partial}^{2}(\phi^{\epsilon}_t\circ\psi^{\epsilon})(Q)\right]
\displaystyle+\epsilon ~\left[
\overline{\partial}^{2}(\phi^{\epsilon}_t\circ\psi^{\epsilon})(Q)~{\partial}(\phi_t\circ\psi)(Q)\right]
\displaystyle+\left[\overline{\partial}^{2}(\phi^{\epsilon}_t\circ\psi^{\epsilon})(Q)\right]^2\\
\\
\Longrightarrow
\vertiii{\left[(\phi^{\epsilon}\circ\psi^{\epsilon})(Q)-\phi_t\left(Q\right)\right]^2
\displaystyle-\epsilon^2~\left[
 {\partial}(\phi\circ\psi)(Q)\right]^2}_{t,n}
\displaystyle\leq \epsilon^{2+1/2}~\overline{e}(t)~\vert Q\vert~\vert H\vert_+
\end{array}
$$
This ends the proof of the estimate (\ref{estimate-tcl}) when $k=2$. The proof of the theorem is completed.\cqfd

\appendix
\section{Appendix: Some technical proofs} \label{appendix-section}

\subsection*{Proof of corollary~\ref{cor-nabla-Pi}}\label{proof-cor-nabla-Pi}

Arguing as in the proof of corollary 4.13 in~\cite{ap-2016} we have
$$
\overline{\nabla}E_{s,t}[Q,H]
=E_{s,t}(Q+H)~H-\int_s^t E_{u,t}(Q+H)
(\phi_u(Q+H)-\phi_u(Q))E_{s,u}(Q)~du
$$ 
Using (\ref{def-gamma}) and (\ref{estimate-nabla-k-phi}) we check the estimates
$$
\vertiii{ \overline{\nabla}E_{s,t}[Q,H]}\leq 
\overline{g}_{\nu}\left(\Vert Q\Vert+\Vert H\Vert,t-s\right)
$$
Arguing as in the proof of theorem~\ref{lem-frechet-derivatives} we check that $E_{s,t}(Q)$ is smooth and the derivatives can be computed using for any $n\geq 1$ the induction
 \begin{align*}
 \nabla^n E_{s,t}(Q)\cdot H~
=&-\int_s^t E_{u,t}(Q)
\left(\nabla^n \phi_u(Q)\cdot H\right)\,E_{s,u}(Q)\,du+n\,
\left[\nabla^{n-1} E_{s,t}(Q)\cdot H\right]\,H\\
& -~n!~\int_s^t\,\sum_{k+l=n-2} \left[\frac{1}{k+1!}\,\nabla^{k+1} E_{u,t}(Q)\cdot H\right]
\left[\frac{1}{(l+1)!}\,\nabla^{l+1} \phi_u(Q)\cdot H\right]E_{s,u}(Q)\,du
\end{align*}
Using the estimates (\ref{estimate-nabla-k-phi}) we check (\ref{nabla-En}).
Also observe that
\begin{eqnarray*}
 \overline{\nabla}\Gamma_{t}[Q,H]&=&\int_0^t~
[E_u(Q+H)-E_u(Q)]^{\prime}~E_u(Q)~du+
\int_0^t~E_u(Q)^{\prime}~[E_u(Q+H)-E_u(Q)]~du\\
&&
 \displaystyle\hskip3cm+\int_0^t~
[E_u(Q+H)-E_u(Q)]^{\prime}~[E_u(Q+H)-E_u(Q)]~du \end{eqnarray*}
This implies that
\begin{eqnarray*}
\frac{1}{n!}\nabla^n\Gamma_t(Q)\cdot H&=&\int_0^t~
\left[\frac{1}{n!}\nabla^n E_u(Q)\cdot H\right]^{\prime}~E_u(Q)~du+
\int_0^t~E_u(Q)^{\prime}~\left[\frac{1}{n!}\nabla^n E_u(Q)\cdot H\right]~du\\
&&\hskip.3cm
 \displaystyle+\int_0^t~\sum_{k+l=n-2}~\left[\frac{1}{(k+1)!}\nabla^{k+1} E_u(Q)\cdot H\right]^{\prime}
~\left[\frac{1}{(l+1)!}\nabla^{l+1} E_u(Q)\cdot H\right]~du  \end{eqnarray*}
The proof of (\ref{estimate-nabla-Pi}) now follows standard computation, thus it is skipped. This ends the proof of corollary~\ref{estimate-nabla-Pi}.
\cqfd

\subsection*{Proof of proposition~\ref{propM-bounds-ref}}\label{propM-bounds-ref-proof}

Firstly, we check (\ref{ref-M-bound-mn}) for $m=1$. We set $(P_t^{(1)},Q_t^{(1)})=(P_t,Q_t)$.
By (\ref{ab-interm-2}) we have
$$
\begin{array}{l}
  \displaystyle dM_t=(P_t\,\bullet \,d\Wa_t\,\bullet Q_t)_{ sym}\\
\\
\Longrightarrow 4\,\tr\left[\partial_t\,\langle\,M\,\vert\,\otimes
\,\vert\, M\rangle^{\,\sharp}_t\right]=2~(\tr(P_tQ_t))^{2}
+2\tr(P_t^2)\tr(Q_t^2)\leq 4 \Vert P_t\Vert_F^{2}\Vert Q_t\Vert_F^{2}
\end{array}
$$
We check (\ref{ref-M-bound-mn}) using the estimates
\begin{eqnarray*}
\vertiii{M}_{t,2n}^2
&\leq &c~\EE\left[\left(\int_0^t~\Vert P_s\Vert^{4}~ds\right)^{n}\right]^{1/(2n)}
~
\EE\left[\left(\int_0^t~\Vert Q_s\Vert^{4}~ds\right)^{n}\right]^{1/(2n)}
\end{eqnarray*}
We let  ${\delta}=(\delta_i)_{i\geq 1}$ be a sequence of independent $\{-1,+1\}$-valued Rademacher random variables.
In this notation, the martingale defined in (\ref{ref-M-bound-mn}) satisfies the polarization formula
$$
M_t=\EE\left[M_t^{\delta}
\vert \Fa_t\right]\Longrightarrow \Vert M_t\Vert\leq \EE\left[\Vert M_t^{\delta}\Vert
\vert \Fa_t\right]
$$
with
$$
\begin{array}{l}
  \displaystyle
M_t^{\delta}:=\int_0^t~\left(\left[\sum_{1\leq i\leq m}\delta_i~P_s^{(i)}\right]\, \,d\Wa_s\, \left[
\sum_{1\leq j\leq m}\delta_j~Q_s^{(j)}\right]\right)_{ sym}
\end{array}
$$
Arguing as above we find that
\begin{eqnarray*}
\vertiii{M}_{t,2n}^2
&\leq& c~ \sum_{1\leq i,j\leq m}\left[\int_0^t~\EE(\Vert P_s^{(i)}\Vert_F^{4n})^{1/n}~ds\right]^{1/2}
~
\left[\int_0^t~\EE(\Vert Q_s^{(j)}\Vert_F^{4n})^{1/n}~ds\right]^{1/2}
\end{eqnarray*}
from which we prove (\ref{ref-M-bound-mn}).
Using (\ref{form-pq-ref}) and (\ref{trQ2}) we also have
$$
\begin{array}{l}
\displaystyle\tr\left[\left((A_s^{\prime}A_s)\otimes (B_sB_s^{\prime})\right)~
\partial_s\langle M^{\epsilon}\,\vert\,\otimes\,\vert \,M^{\epsilon}\,\rangle^{\,\sharp}_s
\right]
\displaystyle\leq \Vert A_s\Vert_F^2~\Vert B_s\Vert_F^2~\sum_{1\leq i,j\leq m}\Vert P_t^{(i)}\Vert_F^{2}~
\Vert Q_t^{(j)}\Vert_F^{2}
\end{array}
$$
This ends the proof of (\ref{ref-M-bound-mn-2-integral}).
Using Holder inequality we check that
$$
\begin{array}{l}
\displaystyle\vertiii{M}_{t,2n}^2
\displaystyle\leq c~\EE\left[\left(\int_0^t~\Vert A_s\Vert^{8}~ds\right)^{n}\right]^{1/(4n)}\EE\left[\left(\int_0^t~\Vert B_s\Vert^{8}~ds\right)^{n}\right]^{1/(4n)}\\
\hskip4cm\displaystyle\sum_{1\leq i,j\leq m}\EE\left[\left(\int_0^t~\Vert P_s^{(i)}\Vert^{8}~ds\right)^{n}\right]^{1/(4n)}
~
\EE\left[\left(\int_0^t~\Vert Q_s^{(j)}\Vert^{8}~ds\right)^{n}\right]^{1/(4n)}\end{array}
$$
The end of the proof of (\ref{ref-M-bound-mn-2}) is now immediate. This ends the proof of the proposition.
\cqfd

\subsection*{Proof of lemma~\ref{lem-pre-cv}}\label{some-ang-prop-proof}

The proof of lemma~\ref{lem-pre-cv} is itself based on the following technical lemma.

\begin{lem}\label{some-ang-prop}
For any matrix valued martingales
$$
M:=P\bullet \Wa\bullet Q\quad \mbox{and}\quad \overline{M}:=\overline{P}\bullet \Wa\bullet\overline{Q}
$$
we have the angle bracket formula
    \begin{equation}\label{ab-interm-2}
      \begin{array}{l}
  \displaystyle 4\,\partial_t\,\langle\,(M-\overline{M})_{ sym}\,\vert\,\otimes
\,\vert\, (M-\overline{M})_{ sym}\rangle^{\,\sharp}_t\\
\\
=(P_tQ_t-\overline{P}_t\,\overline{Q}_t)^{\otimes 2}+(Q_tP_t-\overline{Q}_t\,\overline{P}_t)^{\otimes 2}
+\left[P_t\otimes Q_t-\overline{P}_t\otimes \overline{Q}_t\right]^2+\left[Q_t\otimes P_t-\overline{Q}_t\otimes \overline{P}_t\right]^2\\
\\
\hskip.3cm+2\left\{
(\overline{P}_t\,\overline{\otimes}\, \overline{Q}_t)\left[
\overline{Q}_t\,\otimes\,\overline{P}_t-Q_t\,\otimes\, P_t\right]\right\}_{ sym}+2\left\{
(\overline{Q}_t\,\overline{\otimes}\, \overline{P}_t)\left[
\overline{P}_t\,\otimes\,\overline{Q}_t-P_t\,\otimes\, Q_t\right]\right\}_{ sym}
     \end{array}   
        \end{equation}
\end{lem}

\proof
To simplify the presentation, we drop the time subscript and we write $(P,\overline{P},Q,\overline{Q})$ instead of
$(P_t,\overline{P}_t,Q_t,\overline{Q}_t)$.
Using (\ref{prop-ra-lg-ll-gg-sharp}) and the commutation property (\ref{prod-prop}) we check that
$$
\begin{array}{l}
  \displaystyle\partial_t\,\langle~M\,\vert\,\otimes
\,\vert\, \overline{M}~\rangle^{\,\sharp}=(PQ)\otimes (\overline{P}\,\overline{Q})
\quad\mbox{and}\quad
  \displaystyle\partial_t\,\langle~M\,\vert\,\otimes
\,\vert\, \overline{M}^{\,\prime}~\rangle^{\,\sharp}=(P\overline{P})~\overline{\otimes}~ 
(Q\overline{Q})\\
\\
\Longrightarrow   \displaystyle\partial_t\,\langle~M-\overline{M}\,\vert\,\otimes
\,\vert\, M-\overline{M}~\rangle^{\,\sharp}=(PQ-\overline{P}\,\overline{Q})^{\otimes 2}
\end{array}
$$
In the same vein we check that
      $$
      \begin{array}{l}
  \displaystyle\partial_t\,\langle~M-\overline{M}\,\vert\,\otimes
\,\vert\, \left(M-\overline{M}\right)^{\prime}~\rangle^{\,\sharp}\\
\\
=\left[P\otimes Q-\overline{P}\otimes \overline{Q})\right]^{2}+
(\overline{P}~\overline{\otimes}~ \overline{Q})\left[
(\overline{Q}\otimes\overline{P})-(Q\otimes P)\right]
+\left[(\overline{P}\otimes \overline{Q})-\left(P\otimes Q\right)\right](\overline{P}~\overline{\otimes}~\overline{Q})
    \end{array}
    $$
 The end of the proof (\ref{ab-interm-2}) now follows elementary manipulations. \cqfd

{\bf Proof of lemma~\ref{lem-pre-cv}: }Using (\ref{ab-interm-2}) and the norm estimates (\ref{norm-estimates-ref-frown}) we check that
$$
      \begin{array}{l}
  \displaystyle 2\,\Vert\partial_t\,\langle\,(M-\overline{M})_{ sym}\,\vert\,\otimes
\,\vert\, (M-\overline{M})_{ sym}\rangle^{\,\sharp}\Vert_F\\
\\
  \displaystyle\leq 
 4 ~\tr(\overline{P})~ \Vert Q^{1/2}-\overline{Q}^{1/2} \Vert_F^2+4~\tr(
Q)~
 \Vert P^{1/2}-\overline{P}^{1/2}\Vert_F^2\\
\\\hskip1cm+
2~\sqrt{\tr(\overline{P})~\tr(\overline{Q})}~\sqrt{\tr(Q)}~
\Vert P^{1/2}-\overline{P}^{1/2}\Vert_F
+
2~\tr(\overline{P})~\sqrt{\tr(\overline{Q})}~\Vert Q^{1/2}-\overline{Q}^{1/2} \Vert_F
     \end{array}
     $$
     Using (\ref{square-root-key-estimate}) we check the estimate
     $$
      \begin{array}{l}
  \displaystyle \Vert\partial_t\,\langle\,(M-\overline{M})_{ sym}\,\vert\,\otimes
\,\vert\, (M-\overline{M})_{ sym}\rangle^{\,\sharp}\Vert_F\\
\\
  \displaystyle\leq \left[
 2~\tr(\overline{P})~ \lambda^{-1}_{min}(\overline{Q})~ \left( \Vert Q \Vert_F+ \Vert \overline{Q} \Vert_F\right)
+
\tr(\overline{P})~\sqrt{\tr(\overline{Q})}~\lambda^{-1 / 2}_{min}(\overline{Q})~\right] \Vert Q-\overline{Q} \Vert_F\\
\\\hskip.1cm+\left[
\sqrt{\tr(\overline{P})~\tr(\overline{Q})}~\sqrt{\tr(Q)}~
\lambda^{-1/2}_{min}(\overline{P})
+2~\tr(
Q)~
\lambda^{-1}_{min}(\overline{P})~\left( \Vert P \Vert_F+ \Vert\overline{P} \Vert_F\right)\right]
 \Vert P-\overline{P} \Vert_F
     \end{array}
     $$ 
 The proof of the lemma is now easily completed. \cqfd

\end{document}